\documentclass[journal,twoside,web]{ieeecolor}

\usepackage{generic}
\usepackage{cite}
\usepackage{amsmath,amssymb,amsfonts,bbm} 
\usepackage{mathdots,yhmath,bm}
\usepackage{xfrac}
\usepackage{algorithmic}
\usepackage{graphicx,subfig,epstopdf}
\usepackage{algorithm,algorithmic}
\usepackage{comment}
\usepackage{hyperref}
\hypersetup{hidelinks=true}
\usepackage{textcomp}

\usepackage{mathtools}

\newtheorem{theorem}{Theorem}
\newtheorem{lemma}{Lemma}
\newtheorem{proposition}{Proposition}
\newtheorem{corollary}{Corollary}

\newtheorem{problem}{Problem}
\newtheorem{definition}{Definition}
\newtheorem{assumption}{Assumption}
\newtheorem{remark}{Remark}

\DeclareMathOperator{\tr}{Tr}

\DeclareMathOperator{\diag}{diag}

\DeclareMathOperator{\svec}{svec}
\DeclareMathOperator{\mvec}{vec}

\newcommand{\norm}[1]{\left\lVert#1\right\rVert}

\usepackage{xcolor}

\usepackage{algorithmic}

\newlength\myindent
\begin{document}
\title{Adaptive Inverse Reinforcement Learning with Online Off-Policy Data Collection}
\author{Yibei Li, Yuexin Cao, Zhixin Liu, and Lihua Xie, \IEEEmembership{Fellow, IEEE}
\thanks{The work was supported by Natural Science Foundation of China
under Grant T2293772.}
\thanks{Y. Li and Z. Liu are with State Key Laboratory of Mathematical Sciences, Academy of Mathematics and Systems Science, Chinese Academy of Sciences, China (e-mail: yibeili@amss.ac.cn; lzx@amss.ac.cn). Y. Cao is with  Department of Mathematics, KTH Royal Institute of Technology, Stockholm 10044, Sweden (e-mail: yuexin@kth.se). L. Xie is with School of Electrical and Electronic Engineering, Nanyang Technological University, Singapore (e-mail: elhxie@ntu.edu.sg).  }
}

\maketitle

\begin{abstract}
In this paper, the inverse reinforcement learning (IRL) problem is addressed to reconstruct the unknown cost function underlying an observed optimal policy in a model-free manner, whose online adaptation with completely off-policy system data still remains unclear in the literature. Without prior knowledge of the system model parameters, an adaptive and direct learning rule for the cost parameter is proposed using online off-policy system data, which only needs to satisfy the mild persistently exciting condition in the general data-driven paradigm. The adaptive and online IRL algorithm is achieved by designing full Nesterov-Todd (NT)-step primal-dual interior-point iterations.  Despite solving a nonlinear and time-varying semi-definite program (SDP), the influence of system noise is rigorously analyzed, and the proposed online algorithm is shown to achieve a sublinear convergence. The proposed method is further generalized to nonlinear IRL based on differential dynamic programming. The gradient of the loss function is directly obtained via a backward pass, which eliminates the need to repeatedly solve forward RL problems as in conventional bi-level IRL frameworks. Finally, the efficiency and effectiveness of the proposed algorithms are demonstrated by numerical examples.
\end{abstract}


\section{INTRODUCTION}
\IEEEPARstart{O}{ver} the past few decades, reinforcement learning has emerged as a powerful tool to address complex decision-making problems, driven by its notable efficacy in optimizing a cumulative reward through repeated interaction with a dynamic environment. On the contrary, inverse reinforcement learning aims to infer the underlying cost function that governs the observed optimal policies, which has demonstrated remarkable capabilities across diverse disciplines such as robotics, adversarial games, and human-robot interaction \cite{UCHIBE2021138,KollmitzIROS}.

Existing IRL frameworks either focus on Markov decision processes with various learning techniques (such as MaxEnt \cite{Ziebart2008} and deep neural networks \cite{FinnGuided2016}), or on state-space models that are widely applied in the control community. This paper falls within the latter line of research, which is closely related to the topic of inverse optimal control (IOC). The classic IOC problem, first proposed in \cite{kalman1964linear}, is built on assuming precise knowledge of the system dynamics. For example, well-posedness of the inverse problem is discussed in \cite{li2020continuous,zhang2019inverse,yu2022inverse} and the unknown cost parameters are usually reconstructed within residual minimization frameworks \cite{JinPDP,englert2017inverse}.

In practice, precise system models are rarely known a priori, which are usually estimated from measured data via system identification (SysID). Consequently, a natural routine for model-free IOC/IRL is rooted in an indirect framework, namely, SysID followed by model-based methods, where a least-squares estimate of the system model is first obtained from a batch of offline data and model-based IOC techniques are then applied by using the estimated system model as the ground truth \cite{SELF2022110242,Rushikesh2018ACC}. As has been broadly discussed in the data-driven design paradigm, one major concern of the indirect method lies in the inconsistent optimality of the two sequential subproblems that are usually difficult to integrate \cite{HJALMARSSON2005393}.

\begin{figure}
\begin{center}
\includegraphics[width=8.0cm]{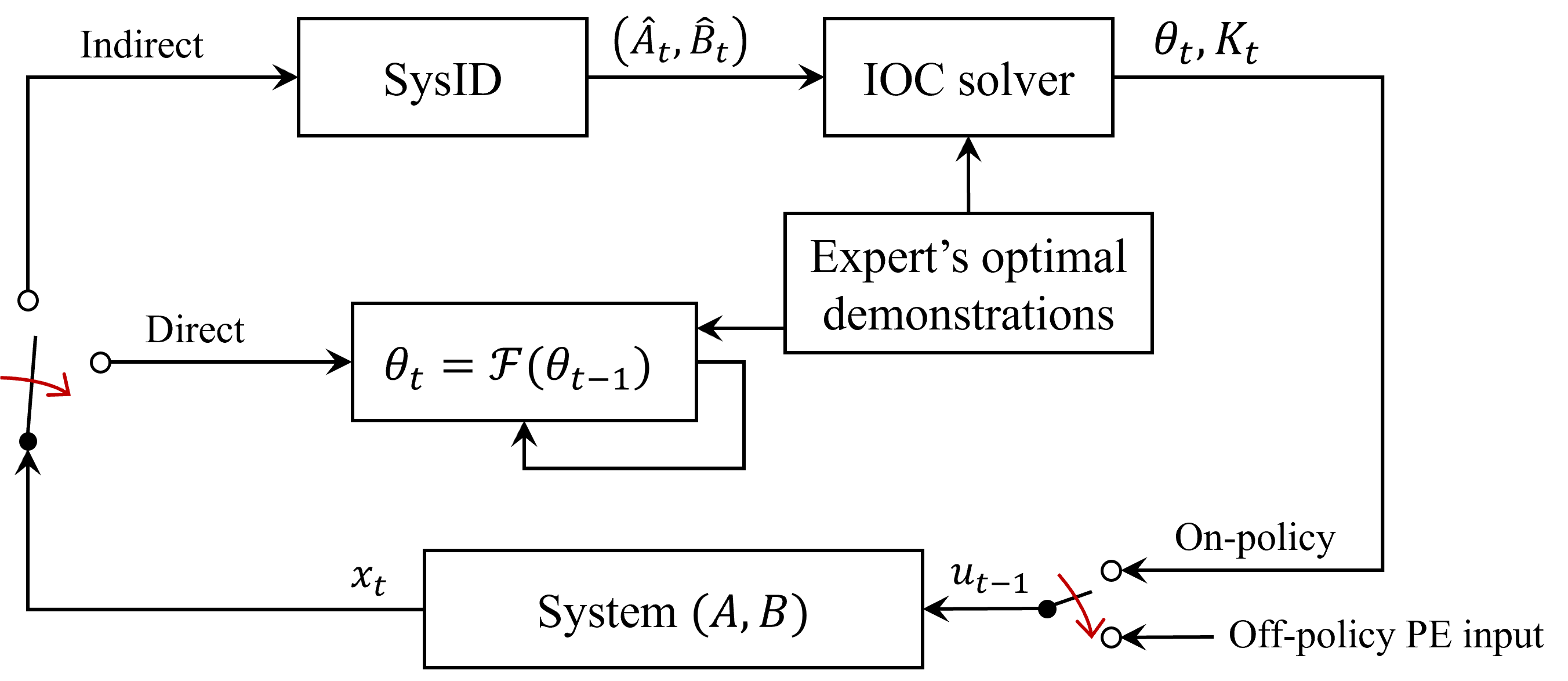} 
  \caption{Illustration of direct/indirect IRL algorithms with on-/off-policy online data, where $\mathcal{F}$ has an explicit form.} 
  \label{flowchart}
\end{center}  
\vspace{-0.7cm}
\end{figure}

Different from the indirect framework, the goal of direct data-driven IRL is to design an end-to-end learning rule of the cost parameter from a single episode of persistently exciting (PE) online data, where its real-time adaptation rule has recently emerged as a critical topic. As illustrated in Fig. \ref{flowchart}, when new data is collected, the indirect method alternates between SysID and solving a complete IOC problem, while the direct method only applies an explicit one-step updating rule of the cost parameter, which is more computationally efficient and conceptually simpler. Although recent studies have achieved direct and adaptive data-driven optimal control \cite{ZhaoData2025}, the issue of adaptive learning still remains largely unexplored in the inverse problem.

Among various types of cost structures, study on the linear quadratic regulator (LQR) has aroused considerable interest, not only because of its analytic tractability, but also lies in its foundational role in many nonlinear optimization schemes. An active line of related research is built on the inverse $Q$-learning framework, where techniques such as policy iterations are applied to derive a data-driven implementation for the model-based IOC conditions. Examples can be found in imitating learning \cite{xue2021inverse}, tracking control \cite{XueIRL2022}, output-feedback differential game \cite{MartirosyanCSL2024} and a series of follow-up research. To identify the cost parameter underlying the observed demonstrations in a model-free manner, the core idea of data-driven IRL is to collect online episodes of system data to parameterize the system dynamics. 
Most of the aforementioned methods require the system data to be collected in an on-policy manner to incorporate the IRL iterations, where the expert needs to persistently apply the current policy estimate to generate online system data. However, in many scenarios, the expert performs independently and cannot be excited arbitrarily in a predefined manner. To alleviate such issue, some attempts design online IRL algorithms that only rely on expert demonstrations \cite{SELF2022110242,asl2022online}, whereas the input dynamics matrix is required as a priori. %

On the other hand, conventional IRL algorithms are usually built on a bi-level framework, where the cost parameter is updated in the outer loop and forward RL is repeatedly solved in the inner loop at the current candidate cost parameter \cite{PERRUSQUIA202233,XueIRL2022}. To alleviate such computational complexity, some end-to-end computations of the cost gradient are derived by directly differentiating the optimality conditions in the inner loop \cite{JinPDP,CaoIRL2025}, where all the unknown parameters in the cost function and system dynamics are updated simultaneously. However, informativity issues may arise in data collection when the system dynamics is completely unavailable as in our problem.

In this paper, we propose a novel direct and adaptive IRL algorithm for LQRs with online collection of off-policy system data. Different from existing online IRL such as \cite{SELF2022110242,asl2022online} that collect online expert demonstrations to continue the learning iterations, we aim to further improve the data efficiency by only using a finite optimal demonstration history, and to design an adaptive rule to update the cost estimate with completely off-policy online system data that only needs to satisfy the basic PE assumption in data-driven paradigms. Based on the differential dynamic programming (DDP) framework, our proposed methods are further generalized to nonlinear cases. Our main contributions are summarized as follows.

\begin{itemize}
    \item A novel online IRL algorithm is designed to adaptively improve the cost parameter estimate by collecting completely off-policy data, where neither the current policy iteration nor the true optimal controller has to be applied to excite the expert system. Adaptive and direct one-step recursions are derived to avoid a complete RL procedure as in conventional bi-level methods.
    \item The influence of system noise on the online iterations is quantified and non-asymptotic convergence results are established, which has been rarely discussed in the existing IRL results. Besides, the choice of starting points and stepsizes does not need any prior knowledge of the true cost parameter as in \cite{xue2021inverse} and similar works. Despite solving a time-varying and nonlinear SDP, the online algorithm converges in a sublinear rate by $\mathcal{O}(1/k)$.
    \item To alleviate the informativity issue in updating both the cost and system parameters simultaneously, our method only involves a cost parameter estimator that in essence converges to the solution of the indirect method. DDP-based IRL is derived for nonlinear costs and its end-to-end computation of the error gradient is proposed, making it possible to analyze online adaptations in future work. 
\end{itemize}

The rest of paper is organized as follows. The inverse data-driven optimal control problem is formulated in Section \ref{sec:prob}. Off-policy online cost learning is studied for the LQR in Section \ref{sec:LQR} and is generalized to the nonlinear case via a DDP-based framework in Section \ref{sec:nonlinear}. Numerical study is provided in Section \ref{sec:simulation} while Section \ref{sec:conclusion} gives some concluding remarks.

\textit{Notations:} $\|\cdot\|$ denotes the $l_2$ norm of a vector or the Frobenius norm of a matrix. Let $\|\cdot\|_{2}$ denote the spectral norm of a matrix and $\langle \cdot,\cdot \rangle$ denote the Frobenius inner product of two matrices. Denote $I_{n}$ as the $n \times n$ identity matrix and $\mathcal{I}_n$ as the identity matrix of order $n(n + 1)/2$. Denote $\otimes$ as the Kronecker product, $\tr(A)$ as the trace of a square matrix $A$ and $\overline{A}:=\frac{1}{2}(A+A^{\top})$.  Let $\diag(X_1,...,X_n)$ denote the block matrix with $X_1,...,X_n$ on its diagonal, where $X_i$ can be either a scalar or a square matrix. For any $x \in \mathbb{R}^n$ and $A \in \mathbb{R}^{n \times n}$, let $\norm{x}_A^2$ be short for $x^{\top}Ax$. We use $\bar{\sigma}(\cdot)$ and $\underline{\sigma}(\cdot)$ for the largest and smallest singular value of a matrix respectively. Let $\mathbb{S}^n$, $\mathbb{S}^n_+$ and $\mathbb{S}^n_{++}$ denote the set of $n \times n$ symmetric matrices, positive semi-definite matrices and positive definite matrices, respectively.  We also use $A \succeq 0$ and $A \succ 0$ to indicate $A \in \mathbb{S}^n_+$ and $A \in \mathbb{S}^n_{++}$, respectively.

\section{Problem Formulation}\label{sec:prob}
In this paper, we consider the discrete-time linear time-invariant system in its general form:
\begin{equation}\label{eq:sys}
    x_{k+1}=Ax_k+Bu_k+w_k,
\end{equation}
where $x_k\in\mathbb{R}^n$ and $u_k\in\mathbb{R}^m$ are the state and control input at time $k$, respectively. Assume that the independent noise sequence $\{w_k\}_{k\geq 0}$ satisfies $\mathbb{E}[w_k]=0$, $\norm{w_k}\leq\kappa$ for some constant $\kappa>0$.

Consider the stochastic linear quadratic problem
\begin{equation}\label{eq:LQR}
\begin{aligned}
    \min_{u}~~&  \mathop{\lim \sup}_{T\rightarrow\infty} \frac{1}{T}\sum_{k=0}^{T-1}(x_k^{\top}Qx_k+\norm{u_{k}}^2) \\
    s.t.~~&x_{k+1}=Ax_k+Bu_k+w_k,\quad x_0~\text{is given},
\end{aligned}
\end{equation}
Let $P\succ0$ uniquely solve the algebraic Riccati equation (ARE):
\begin{equation}\label{eq:are}
  \begin{aligned}
    P = A^TPA+Q-A^TPB(B^TPB+I_{m})^{-1}B^TPA.
\end{aligned}   
\end{equation}    
Then the optimal control to \eqref{eq:LQR} is given by $u_k^*=-\bar{K}x_k^*$ ~\cite{ChenBook} with 
\begin{equation}
    \bar{K}=(B^{\top}PB+I_{m})^{-1}B^{\top}PA.
\end{equation}

As will be studied in Section \ref{sec:nonlinear}, the goal of general nonlinear IRL is to infer the unknown cost parameters from optimal demonstrations of $\{x_k^*,u_k^*\}_{k=0}^{T-1}$. Since in LQRs the optimal feedback matrix is reduced to a constant matrix that can be uniquely pre-estimated from a finite dataset of $(x_k^*,u_k^*)$, the inverse LQR (ILQR) problem usually aims to reconstruct the penalty matrix $Q$ directly from $\bar K$. However, without prior knowledge of the model parameter $(A,B)$, data-driven IRL relies on collecting online system data  (that may not be optimal) to characterize the input-output behavior of the linear system and to develop an adaptive rule for cost estimation.

Given any system trajectories on a time interval $[0,t]$, we define the following notations:
\begin{align*}
    & X_{0,t}:=[x_0~x_1~\cdots~x_{t-1}],~~U_{0,t}:=[u_0~u_1~\cdots~u_{t-1}],\\
    & X_{1,t}:=[x_1~x_2~\cdots~x_{t}],~~W_{0,t}:=[w_0~w_1~\cdots~w_{t-1}].
\end{align*}

By the Willems fundamental lemma, it is well-known that the input-output data provides a data-driven parametrization of a linear system under the persistently exciting (PE) condition.
\begin{definition}[PE condition]
    An input-state series with length $t$ is called persistently exciting if the data matrix
    \begin{equation}\label{eq:PE_Z}
        Z_{0,t}:=\begin{bmatrix}
            X_{0,t}\\U_{0,t}
        \end{bmatrix}\in\mathbb{R}^{(n+m)\times t},
    \end{equation}
    has full row rank.
\end{definition}

The PE rank condition is commonly used in data-driven approaches for linear systems. For stochastic systems, various notions have been developed to quantify the PE level, among which we adopt the one used in \cite{ZhaoData2025,coulson2022quantitative} as follows, which guarantees the rank condition in \eqref{eq:PE_Z}.

\begin{assumption}\label{assump:PE}
Let $t_0$ denote the length of the initial dataset. For any $t\geq t_0$ along the online state trajectory, there exists constant $\gamma>0$ such that $\underline{\sigma}(\mathcal{H}_{n+1}(U_{0,t}))\geq\gamma\sqrt{t(n+1)}$, where $\mathcal{H}_{n+1}(U_{0,t})$ is the corresponding Hankel matrix.
\end{assumption}

In the sequel, we refer to $\gamma/\kappa$ as the signal-to-noise ratio (SNR).
Using PE data to approximate the system model, we aim to design a direct and adaptive data-driven algorithm to solve the ILQR problem with online off-policy data as follows, which will be generalized to nonlinear IRL in Section \ref{sec:nonlinear}.
\begin{problem}\label{prob:IOC} 
Given the optimal feedback control matrix $\bar{K}$ that solves \eqref{eq:LQR}, the goal is to recover the underlying cost matrix $Q$ without prior knowledge on the system matrices $(A,B)$. Given an initial system dataset $(X_{0,t_0},U_{0,t_0},X_{1,t_0})$, a new pair of online off-policy system data $(x_t,u_t)$ is collected at each time instant $t> t_0$, based on which we aim to develop a recursive learning rule for the cost matrix.
\end{problem}
\begin{remark}
    Unlike existing on-policy methods that require system data to be collected in a specific manner to achieve model-free implementation of the IRL algorithm, the proposed method significantly improves data efficiency by addressing the fundamental issue that how arbitrary PE data (even collected in other tasks) can be leveraged in the adaptive IRL. In our results, without exciting the expert system with the policy estimate in each iteration, the online off-policy data $(x_t,u_t)$ only needs to satisfy the mild PE condition under Assumption \ref{assump:PE}, which is much weaker and more reasonable. 
\end{remark}

\section{Adaptive Off-Policy Data-Driven ILQR}\label{sec:LQR}
In this section, an adaptive algorithm to reconstruct the unknown LQ cost function is proposed with online off-policy data. A data-driven formulation of the ILQR problem is derived in Section \ref{sec:iLQR_modelfree}, based on which an offline algorithm is first proposed in Sections \ref{sec:alg_offline} \& \ref{sec:SDP_convergence} and further developed into an online adaptive algorithm with guaranteed convergence in Section \ref{sec:online}.
\subsection{Data-Driven Formulation of ILQR}\label{sec:iLQR_modelfree}
In this part we propose a data-driven formulation of the ILQR problem in the deterministic case. The noise influence will be further studied in Section \ref{sec:online}. Define
\begin{equation}\label{eq:HQ}
\setlength{\arraycolsep}{3pt}
H(Q)=\begin{bmatrix*}[c] Q+A^{\top}PA & A^{\top}PB \\
        B^{\top}PA & I_{m}+B^{\top}PB
\end{bmatrix*} \!:=\!\begin{bmatrix}
H_{xx} & H_{ux}^{\top} \\
H_{ux} & H_{uu}
\end{bmatrix},
\end{equation}
where $P$ satisfies the ARE in \eqref{eq:are} with parameter $Q$. It is well-known that the optimal value function for the forward LQR is given by $V(x_k)=x_k^{\top}Px_k$. Define the $\mathcal{Q}$-function as 
\begin{equation}
\begin{aligned}
    &\mathcal{Q}(x_k,u_k)=x_k^{\top}Qx_k+u_k^{\top}u_k+V(x_{k+1})\\
    =&x_k^{\top}Qx_k+u_k^{\top}u_k+(Ax_k+Bu_k)^{\top}P(Ax_k+Bu_k)\\
    =&\begin{bmatrix}
        x_k^{\top} & u_k^{\top}
    \end{bmatrix}H\begin{bmatrix}
        x_k \\ u_k
    \end{bmatrix}.
\end{aligned}
\end{equation}
Then the optimal feedback control indeed solves $u_k^*=\mathop{\arg\min}_{u_k} \mathcal{Q}(x_k,u_k)=-\bar{K}x_k=-H_{uu}^{-1}H_{ux}x_k$. 

It is noteworthy that when solving the inverse problem, the penalty matrix should be constrained within a feasible set such that the corresponding forward problem is well-defined, namely, there always exists a unique optimal and stabilizing control for any $x_0$, which is involved implicitly in the following proposition. 

\begin{proposition}\label{prop:ARE_H}
    Given $(A,B)$ and $Q\succeq 0$, if there exists $H\succ 0$ that solves
    \begin{equation}\label{eq:ARE_H}
    H=\begin{bmatrix}
        Q & 0 \\0 & I_m
        \end{bmatrix}+\Phi(H)^{\top}H\Phi(H), 
    \end{equation}
    and $K:={H}_{uu}^{-1}{H}_{ux}$ such that $A-BK$ is stable, where $H$ is partitioned by
    \begin{equation}\label{eq:H_struct}
    H=\begin{bmatrix}
    H_{xx} & H_{ux}^{\top} \\
    H_{ux} & H_{uu}
    \end{bmatrix},
    \end{equation}
    and
    \begin{equation*}
    \Phi(H):=\begin{bmatrix}
        A & B \\-H_{uu}^{-1}H_{ux}A & -H_{uu}^{-1}H_{ux}B
    \end{bmatrix},
    \end{equation*}
    then the forward LQR with penalty matrix $Q$ is well-defined. Furthermore, $H$ is uniquely given by $H(Q)$ as defined in \eqref{eq:HQ}.
\end{proposition}

\proof
    Define $c(x_k, u_k) = x_k^{\top}Qx_k+u_k^{\top}u_k$ and $V(x_k)=[x_k^{\top}~(-Kx_k)^{\top}]{H}[x_k^{\top}~(-Kx_k)^{\top}]^{\top}$. Note that $H\succ 0$ implies $H_{uu}\succ 0$. Following similar ideas as in \cite{lopez2023efficient}, by substituting $x_{k+1}=Ax_k+Bu_k$ and \eqref{eq:ARE_H}, one can first show that 
    \begin{equation*}
    \begin{aligned}
        \mu_k&=\mathop{\arg\min}_{u_k} \lbrace c(x_k, u_k)+ V(x_{k+1}) \rbrace \\
        &=\mathop{\arg\min}_{u_k} \begin{bmatrix}
            x_k^{\top} & u_k^{\top}
        \end{bmatrix}{H}\begin{bmatrix}
            x_k \\ u_k
        \end{bmatrix}
        =-{H}_{uu}^{-1}{H}_{ux}x_k=-Kx_k,
    \end{aligned}
    \end{equation*}
    and $V(x_k)$ satisfies the Bellman equation:
    \begin{equation}\label{eq:HJBE}
    \begin{aligned}
       & \min_{u_k} \lbrace c_k(x_k, u_k)+ V(x_{k+1}) \rbrace\\= &[x_k^{\top}~(-Kx_k)^{\top}]{H}[x_k^{\top}~(-Kx_k)^{\top}]^{\top} = V(x_k).
    \end{aligned}
    \end{equation}
Let $P:=[I,-K^{\top}]{H}[I,-K^{\top}]^{\top}=H_{xx}-H_{ux}^{\top}H_{uu}^{-1}H_{ux}$.
Then $V(x_{k})=x_k^{\top}Px_{k}$ and \eqref{eq:HJBE} indicate that $P$ satisfies the ARE with parameter $Q$.
By the property of Schur complement, $H\succ 0$ further implies 
$P=H_{xx}-H_{ux}^{\top}H_{uu}^{-1}H_{ux}\succ 0$.
Hence, $V(x_k)$ is positive definite and radially unbounded. Furthermore,  since the strong solution $P$ is both stabilizing and positive definite, it is indeed the unique positive-definite solution of the ARE. Therefore, the forward LQR with penalty matrix $Q$ admits a unique optimal control $\mu_k=-Kx_k$ for any $x_0$, where $K=(B^{\top}PB+R)^{-1}B^{\top}P$ and $V(x_k)$ is the optimal cost-to-go function. On the other hand, we define
\begin{equation}\label{eq:data_matrix_trueAB}
    \Phi_K:=\begin{bmatrix}
        A & B \\-KA & -KB
    \end{bmatrix},
\end{equation}
    which is stable as $A-BK$ is stable~\cite[Lemma 2]{lopez2023efficient}. Hence, the discrete-time Lyapunov equation $H=\diag (Q,I_{m})+\Phi_K^{\top}H\Phi_K$ admits a unique solution, which is equal to  \eqref{eq:HQ}. 
\endproof
\begin{remark}
For the LQ structure with penalty matrix $Q\succeq 0$, it is well-known that the necessary and sufficient condition for the corresponding ARE to admit a unique positive definite and stabilizing solution is that $(A,B)$ is stabilizable and $(Q, A)$ has no unobservable modes inside, or on the unit circle~\cite{wimmer1996set}. However, direct data-driven implementation of the above constraints in the inverse problem is not easy, which is instead naturally guaranteed by the condition $H\succ 0$ and that ${H}_{uu}^{-1}{H}_{ux}$ is stabilizing in our formulation.
\end{remark}

By Proposition \ref{prop:ARE_H}, an alternative formulation for the ILQR problem based on the $\mathcal{Q}$-function is proposed in the following corollary.
\begin{corollary}\label{cor:ILQR_cond}
    Let $\bar K$ be the optimal feedback matrix of the expert demonstrations. Then $Q\succeq 0$ solves the ILQR problem if and only if there exists $H\in\mathbb{S}_{++}^{m+n}$ such that 
    \begin{subequations} 
        \begin{align}
            & H-\diag(Q,I_{m})-\Phi_{\bar K}^{\top}H\Phi_{\bar K}=0, \label{eq:H_lyp}\\
            & H_{ux}=H_{uu}\bar{K},
        \end{align}
    \end{subequations}
    where $H$ is partitioned as in \eqref{eq:H_struct} and $\Phi_{\bar K}$ is defined in the same way as in \eqref{eq:data_matrix_trueAB} by replacing $K$ with $\bar{K}$.
\end{corollary}

By defining
\begin{equation*}
    \Psi(H)=\begin{bmatrix}
        0_{n\times m
        } \\ I_m
    \end{bmatrix}\begin{bmatrix}
        0_{m\times n} & I_m
    \end{bmatrix}H\begin{bmatrix}
        I_n \\ -\bar{K}
    \end{bmatrix}\begin{bmatrix}
        I_n & -\bar{K}^{\top}
    \end{bmatrix},
\end{equation*}
we can obtain $\norm{H_{ux}-H_{uu}\bar{K}}^2=\langle H,\Psi(H) \rangle$. Then Corollary \ref{cor:ILQR_cond} implies that in the model-based case, the cost matrix $Q$ can be obtained as the global minimum to the following problem with a conic constraint:
\begin{equation} \label{eq:LS_Q_K_given}
 \begin{aligned}
   \mathop {\min }\limits_{Q\in \mathbb{S}^n} {\kern 4pt} & f(Q)=\frac{1}{2}	\langle H(Q),\Psi(H(Q)) \rangle \\
   s.t.{\kern 5pt}  &Q\succeq 0,
 \end{aligned}
\end{equation}
where $H(Q)$ uniquely solves \eqref{eq:H_lyp} at candidate $Q$, and $f(Q)$ will achieve its minimum $f(Q^*)=0$ if and only if $Q^*$ solves the ILQR problem.
\begin{remark}
While solving \eqref{eq:LS_Q_K_given}, the constraint $H\succ 0$ will be naturally satisfied in the primal-dual interior point iterates as shown in Section \ref{sec:alg_offline}. Otherwise, one can also relax $Q\succeq 0$ by $Q \succeq \epsilon I_n$ for some small $\epsilon>0$, which will guarantee $H \succ 0$ as well. However, unlike most existing IRL frameworks built on $Q$-learning, we will show that the proposed formulation enables online iterations with completely off-policy data as well as quantitative analysis of the stochastic system.
\end{remark}

Next, to obtain a scalable data-driven formulation of the cost parameter estimator in \eqref{eq:LS_Q_K_given}, the sample covariance of the input-state data is considered:
\begin{equation}
    \Lambda_t:=\frac{1}{t}Z_{0,t}Z_{0,t}^{\top}.
\end{equation}
Accordingly, the system data matrices are transformed by:
\begin{align*}
    &\bar{X}_{0,t}=\frac{1}{t}X_{0,t}Z_{0,t}^{\top},~\bar{X}_{1,t}=\frac{1}{t}X_{1,t}Z_{0,t}^{\top},\\
    &\bar{U}_{0,t}=\frac{1}{t}U_{0,t}Z_{0,t}^{\top},~\bar{W}_{0,t}=\frac{1}{t}W_{0,t}Z_{0,t}^{\top}.
\end{align*}
\begin{remark}
    Under the PE condition, it is clear that $\Lambda_t$ is always nonsingular. Compared with the original data, incorporating the transformed data has two advantages. Firstly, the dimensions of $\Lambda_t,\bar{X}_{0,t},\bar{X}_{1,t},\bar{U}_{0,t},\bar{W}_{0,t}$ are all invariant with $t$, which guarantees the scalability of the online algorithm. Secondly, defining $z_t=[x_t^{\top}~u_t^{\top}]^{\top}$, the aforementioned quantities can be computed in a recursive and memory-efficient manner \cite{ZhaoData2025} as
\begin{equation}\label{eq:data_recur}
    \begin{aligned}
        &\bar{X}_{1,t+1}=\frac{t}{t+1}\bar{X}_{1,t}+\frac{1}{t+1}x_{t+1}z_t^{\top},\\
        &\Lambda_{t+1}^{-1}=\frac{t+1}{t}(\Lambda_{t}^{-1}-\frac{\Lambda_{t}^{-1}z_tz_t^{\top}\Lambda_{t}^{-1}}{t+z_t^{\top}\Lambda_{t}^{-1}z_t}).
    \end{aligned}
\end{equation}
\end{remark}

Under the PE condition, the least-squares estimate of system parameters in \eqref{eq:sys} is uniquely obtained as $[\hat{A}_t~ \hat{B}_t]=X_{1,t}Z_{0,t}^{\dagger}=\bar{X}_{1,t}\Lambda_t^{-1}$. Therefore, given an offline dataset $(X_{0,t},U_{0,t},X_{1,t})$, the data-driven formulation for \eqref{eq:LS_Q_K_given} yields 
\begin{equation} \label{eq:indirect_IOC}
\begin{aligned}
\setlength{\arraycolsep}{3pt}
   \mathop {\min }\limits_{Q\in \mathbb{S}^n} {\kern 2pt} & f_t(Q):=\frac{1}{2}	\langle H,\Psi(H) \rangle \\
   s.t.{\kern 5pt}  &H\!=\!\begin{bmatrix}
        Q & 0\\0 & I_{m}
        \end{bmatrix}\!+\!\begin{bmatrix}
        \hat{A}_t & \hat{B}_t \\-\bar{K}\hat{A}_t & -\bar{K}\hat{B}_t
    \end{bmatrix}^{\top}\!\!H\!\begin{bmatrix}
        \hat{A}_t & \hat{B}_t \\-\bar{K}\hat{A}_t & -\bar{K}\hat{B}_t
    \end{bmatrix}\!,  \\ 
    &[\hat{A}_t~\hat{B}_t]=\mathop{\arg\min}_{A, B} \norm{X_{1,t}-[A~B]Z_{0,t}}^2,   \\
   & Q\succeq 0, 
\end{aligned}
\end{equation}
which is equivalent to
\begin{equation} \label{eq:Jt}
 \begin{aligned}
\mathop {\min }\limits_{Q\in \mathbb{S}^{n}} {\kern 4pt} & f_t(Q):= \frac{1}{2} \langle H,\Psi(H) \rangle \\
   s.t.{\kern 5pt}  & H=\diag(Q,I_{m})+\hat{\Phi}_t^{\top}H\hat{\Phi}_t, \\
   & Q \succeq 0,
 \end{aligned}
\end{equation}
with $\hat{\Phi}_t$ defined by
    $\hat{\Phi}_t:=\begin{bmatrix}
            I ~ -\bar{K}^{\top}
        \end{bmatrix}^{\top}\bar{X}_{1,t}\Lambda_t^{-1}$.
Under mild assumptions of the SNR, the equality constraint in \eqref{eq:Jt} is still well-defined since $\hat{\Phi}_t$ is stable by Lemma \ref{lem:Phi_t_stable}, whose proof is given in Appendix \ref{append:Phi_t_stable}.
\begin{lemma}\label{lem:Phi_t_stable}
    Assume that the SNR satisfies    
    \begin{equation}\label{eq:SNR_bnd}
        \frac{\gamma}{\kappa}\geq \frac{2}{\zeta} \max\lbrace \omega, \hat{\sigma} \rbrace
    \end{equation}
    with 
    $\hat{\sigma}=(\norm{\bar{K}}_2+1)(\sqrt{\bar{\sigma}(\Sigma_K-I_n)\bar{\sigma}(\Sigma_K)}+\bar{\sigma}(\Sigma_K)),$
    where $\zeta$ depending on $(A,B)$ characterizes the controllability of the system (i.e., as $\rho$ in \cite{coulson2022quantitative}), $\omega$ depending on $A$ characterizes the noise propagation properties (i.e., as $\gamma$ in \cite{coulson2022quantitative}), and $\Sigma_K=\sum_{i=0}^{\infty}(({A}-{B}\bar{K})^{\top})^{i}({A}-{B}\bar{K})^i$.
    Then $\hat{\Phi}_t$ is stable.
\end{lemma}
\begin{remark}
With the offline data $(X_{0,t},U_{0,t},X_{1,t})$, the equivalence of \eqref{eq:indirect_IOC} and \eqref{eq:Jt} implies that the direct data-driven IRL \eqref{eq:Jt} can achieve at least the same accuracy as that of the conventional indirect method,
which in essence first obtains the least-squares estimate of $(A,B)$ via offline data and regards it as the ground truth in the IOC process. However, how to derive an adaptive rule for \eqref{eq:indirect_IOC} with online data is unclear, which will be addressed via \eqref{eq:Jt} in Section \ref{sec:online}.
\end{remark}

\subsection{Primal-Dual Interior-Point Method with Offline Data}\label{sec:alg_offline}
In this part, given the offline data $(X_{0,t},U_{0,t},X_{1,t})$, the direct data-driven ILQR formulation in \eqref{eq:Jt} is solved, which will serve as the foundation for the online adaptive algorithm in the subsequent subsections. To address the conic constraint, a full-step primal-dual interior-point method is developed.

The gradient of the convex loss function $f_t(Q)$ is first derived in Proposition \ref{prop:convF_grad}, based on which the dual problem is obtained in Proposition \ref{prop:dual}. Their proofs are given in Appendices \ref{append:convF_grad} and \ref{sec:append_prop_dual} respectively.
\begin{proposition}\label{prop:convF_grad}
    $f_t(Q)$ is convex in $Q$ and its gradient $\nabla f_t(Q)$ is given by
    \begin{equation}\label{eq:Delta_f}
        \nabla f_t(Q) = \begin{bmatrix}
        I_n & 0
    \end{bmatrix}Z\begin{bmatrix}
        I_n \\ 0
    \end{bmatrix},
    \end{equation}
    where $Z\in \mathbb{S}^{n+m}$ solves
    \begin{equation}\label{eq:Z_syl}
        Z=\overline{\Psi(H(Q))} + \hat{\Phi}_t Z\hat{\Phi}_t^{\top}, 
    \end{equation}
    and $H(Q)$ is defined as the unique solution to the following discrete-time Lyapunov equation at candidate $Q$:
    \begin{equation}\label{eq:H_lyap_t}        H=\diag(Q,I_{m})+\hat{\Phi}_t^{\top}H\hat{\Phi}_t.
    \end{equation}
    \end{proposition}

\begin{proposition}\label{prop:dual}
    Let $S\in\mathbb{S}^n_+$ denote the Lagrange multiplier of \eqref{eq:Jt}, then its dual problem is given as follows.
\begin{subequations} \label{eq:dual_prob}
 \begin{align}
\mathop {\max }\limits_{S,Q,Z,H} {\kern 4pt} & \frac{1}{2}\tr (\begin{bmatrix}-Q & 0 \\0 & I_{m}\end{bmatrix}Z) \notag \\
   s.t.{\kern 5pt}  & S=\begin{bmatrix}
        I_n & 0
    \end{bmatrix}Z\begin{bmatrix}
        I_n \\ 0
    \end{bmatrix},\label{subeq:S_eqn}\\
   &Z=\overline{\Psi(H)} + \hat{\Phi}_t Z\hat{\Phi}_t^{\top}, \label{subeq:Z_eqn}\\
   &H=\diag(Q, I_{m})+\hat{\Phi}_t^{\top}H\hat{\Phi}_t, \label{subeq:H_eqn}\\
   &S\succeq 0.
 \end{align}
\end{subequations}
If $Q$ and $S$ are feasible to both the primal and dual problems, then $\langle Q,S \rangle$ gives the duality gap.
\end{proposition}

\begin{assumption}[Slater regularity condition]\label{assump:regularity}
    Suppose that there exist $Q\succ 0, S\succ 0, H\in\mathbb{S}^{n+m}$ and $Z\in \mathbb{S}^{n+m}$ in the interior of both the feasible sets of the primal and dual problems, namely, satisfying \eqref{subeq:S_eqn}, \eqref{subeq:Z_eqn} and \eqref{subeq:H_eqn}. 
\end{assumption}

In the sequel, we refer to the conditions on $(Q,S,H,Z)$ in Assumption \ref{assump:regularity} as the interior-point condition (IPC) at $t$.
Then the strong duality holds and computing the optimal solutions to the primal problem \eqref{eq:Jt} and dual problem \eqref{eq:dual_prob} is equivalent to solving the following Karush–Kuhn–Tucker (KKT) condition:
\begin{equation}\label{eq:KKT}
        \eqref{subeq:S_eqn}, \, \eqref{subeq:Z_eqn}, \, \eqref{subeq:H_eqn}, \, Q\succeq 0, \, S\succeq 0, \, \langle S,Q \rangle=0. 
\end{equation}

The basic idea of the primal-dual interior point method is to solve a sequence of perturbed KKT conditions (also known as the central path) as follows with $\mu >0$:
\begin{equation}\label{eq:cen_path}
\eqref{subeq:S_eqn}, \, \eqref{subeq:Z_eqn}, \, \eqref{subeq:H_eqn}, \, Q\succ 0, \, S\succ 0, \, SQ = \mu I_n,
\end{equation}
which in fact is the KKT condition to the minimization problem of the following barrier function
\begin{equation}\label{eq:f_mu}
    \mathop {\min }\limits_{Q \in \mathbb{S}^{n}_{++}}~f_{\mu}(Q)=f_t(Q)-\mu\log\det(Q).
\end{equation}

By designing a recursion rule along the central path, the optimal solution to both the primal and dual problems can be obtained as $\mu\rightarrow 0$. The well-posedness of the central path \eqref{eq:cen_path} is first analyzed.
\begin{proposition}
    Under Assumption \ref{assump:regularity}, for any $\mu>0$, there exists a unique solution of $(Q,S,H,Z)\in\mathbb{S}^n_{++}\times\mathbb{S}^n_{++}\times\mathbb{S}^{n+m}\times\mathbb{S}^{n+m}$ to the central path \eqref{eq:cen_path}.
\end{proposition}
\proof
    Firstly, showing the existence of the central path is equivalent to showing that $f_{\mu}(Q)$ admits a minimum in $\mathbb{S}^n_{++}$. Similar to other variants of SDP, the main difficulty lies in that $\mathbb{S}^n_{++}$ is not compact and $f_{\mu}(Q)$ is undefined on its boundary. Under Assumption \ref{assump:regularity}, suppose $\bar{Q}\succ 0, \bar{S}\succ 0, \bar{H}\in\mathbb{S}^{n+m}$ and $\bar{Z}\in \mathbb{S}^{n+m}$ satisfy \eqref{subeq:S_eqn}, \eqref{subeq:Z_eqn} and \eqref{subeq:H_eqn}. Thus $\bar Q$ is also feasible to \eqref{eq:Jt}. Define the level set 
        $\Omega_\mu=\lbrace Q\succ 0: f_{\mu}(Q) \leq f_{\mu}(\bar{Q}) \rbrace$.
    We then show that $\Omega_\mu$ is compact. Since $\bar Z$ and $\bar H$ satisfy \eqref{subeq:Z_eqn}, for any $Q\in\Omega_\mu$ and the corresponding $H$ that satisfies \eqref{subeq:H_eqn}, one can derive that
    \begin{equation}\label{eq:HBarpsi_eqn}
        \begin{aligned}
            &\langle {H}, \overline{\Psi(\bar H)} \rangle = \langle {H}, \bar Z \rangle-\langle {H}, \hat{\Phi}_t \bar Z\hat{\Phi}_t^{\top} \rangle\\
            =&\langle {H}, \bar Z \rangle-\langle \hat{\Phi}_t^{\top}{H}\hat{\Phi}_t, \bar Z \rangle=\langle \diag(Q, I_{m}), \bar Z \rangle.
        \end{aligned}
    \end{equation}
    Hence, $f_{\mu}(Q) \leq f_{\mu}(\bar{Q}) $ implies that
    \begin{equation}\label{eq:Fmu_ineqn}
        \begin{aligned}
            &\frac{1}{2}	\langle H,\Psi(H) \rangle-\langle {H}, \overline{\Psi(\bar H)} \rangle +\langle \diag(Q, I_{m}), \bar Z \rangle\\
            &-\mu\log\det(Q)\leq f_{\mu}(\bar{Q})
        \end{aligned}
    \end{equation}
    Denote $\sigma_s>0$ as the minimal eigenvalue of $\bar S$. Since $Q\succ 0$, it holds that $\langle \diag(Q, I_{m}), \bar Z \rangle \geq \sigma_s\tr(Q)+\tr(\diag(0, I_{m})\bar{Z})$. In addition, we also have that
    \begin{equation*}
        \begin{aligned}
            &\frac{1}{2} \langle H,\Psi(H) \rangle-\langle {H}, \overline{\Psi(\bar H)} \rangle\\
            =&\frac{1}{2} \langle H-\bar{H},\Psi(H-\bar{H}) \rangle - \frac{1}{2}	\langle \bar{H},\Psi(\bar{H}) \rangle \geq - \frac{1}{2}	\langle \bar{H},\Psi(\bar{H}) \rangle.
        \end{aligned}
    \end{equation*}
    Denote $\lambda_1,\cdots,\lambda_n>0$ as the eigenvalues of $Q$. Substituting the above inequalities into \eqref{eq:Fmu_ineqn} then yields
    \begin{equation}\label{eq:eigQ_bound}
    \begin{aligned}
        &f_{\mu}(\bar{Q})+\frac{1}{2}	\langle \bar{H},\Psi(\bar{H}) \rangle-\tr(\diag(0, I_{m})\bar{Z})\\
        \geq & \sigma_s\tr(Q)-\mu\log\det(Q)  
        =\sum_{i=1}^n (\sigma_s\lambda_i-\mu\log\lambda_i).
    \end{aligned}
    \end{equation}
    Note that $\alpha(t):=\sigma_s t-\mu\log t$ is continuous on $(0,\infty)$, which tends to infinity when $t\rightarrow0$ or $t\rightarrow \infty$. Since the left-hand side of \eqref{eq:eigQ_bound} is known and bounded, $\{\lambda_i\}_{i=1}^n$  must lie in a compact set. Thus, $\Omega_\mu$ is compact. Since $f_{\mu}(Q)$ is continuous on $\mathbb{S}^n_{++}$, it must admit a minimum on $\Omega_\mu$. Furthermore, Proposition \ref{prop:convF_grad} implies $f_{\mu}(Q)$ is strictly convex in $Q$ as $\log\det(Q)$ is strictly concave. Therefore, the existence of a unique minimizer to \eqref{eq:f_mu} is guaranteed, which implies that its KKT conditions \eqref{eq:cen_path} admit a unique solution. 
\endproof

To achieve feasible updates following the central path, applying Newton's method to \eqref{eq:cen_path} results in the search direction that is defined by
\begin{subequations}\label{eq:newton}
    \begin{align}
        &\Delta S=\begin{bmatrix}
        I_n & 0
    \end{bmatrix}\Delta Z\begin{bmatrix}
        I_n \\ 0
    \end{bmatrix}, \label{eq:newton_s}\\
   &\Delta Z=\overline{\Psi(\Delta H)} + \hat{\Phi}_t \Delta Z\hat{\Phi}_t^{\top}, \label{eq:newton_z}\\
   &\Delta H=\diag(\Delta Q, 0)+\hat{\Phi}_t^{\top}\Delta H\hat{\Phi}_t, \label{eq:newton_H}\\
   &\Delta QS+Q\Delta S=\mu I_{n}-QS. \label{eq:newton_Q}
    \end{align}
\end{subequations}
To preserve the symmetry of updates, in this paper we apply the Nesterov-Todd (NT) symmetrization scheme \cite{nesterov1998primal} to symmetrize \eqref{eq:newton_Q}. For any $Q\succ 0$ and $S\succ 0$, there exists a unique weighting point $W\succ 0$ such that $W^{-1}QW^{-1}=S$, where $W$ is given by
\begin{equation*}
    W=Q^{\frac{1}{2}}(Q^{\frac{1}{2}}SQ^{\frac{1}{2}})^{-\frac{1}{2}}Q^{\frac{1}{2}}=S^{-\frac{1}{2}}(S^{\frac{1}{2}}QS^{\frac{1}{2}})^{\frac{1}{2}}S^{-\frac{1}{2}},
\end{equation*}
and the spectral geometric mean of $S$ and $Q$ is given by
\begin{equation}\label{eq:def_V}
    W^{-\frac{1}{2}}QW^{-\frac{1}{2}}=W^{\frac{1}{2}}SW^{\frac{1}{2}}:=\sqrt{\mu}V.
\end{equation}
Then \eqref{eq:newton_Q} is equivalent to
\begin{equation}\label{eq:NT_QS}
    W^{-1}{\Delta Q}W^{-1}+\Delta S =\mu Q^{-1}-S,
\end{equation}
where symmetric solutions of $\Delta Q$ and $\Delta S$ can be guaranteed consequently.

Next, the uniqueness of the symmetric search direction is shown in Proposition \ref{prop:NT_direction} by introducing the symmetric Kronecker product and the $\svec$ operator. More details of related definitions and properties are given in Appendix \ref{appendix:svex}. The proof of Proposition \ref{prop:NT_direction} is given in Appendix \ref{append:NT_direction}.

\begin{proposition}\label{prop:NT_direction}
For any $Q,S\succ 0$, there exists a unique search direction $(\Delta S, \Delta Z,\Delta H,\Delta Q)\in\mathbb{S}^n\times \mathbb{S}^{m+n}\times \mathbb{S}^{m+n}\times \mathbb{S}^n$ that satisfies (\ref{eq:newton}a-\ref{eq:newton}c) and \eqref{eq:NT_QS}, which can be uniquely solved by
\begin{equation}\label{eq:newton_svec}
        \begin{bmatrix}
            \mathcal{I}_n & -E & & \\
             & G & F &\\
             & & G^{\top}  & E^{\top}\\
             \mathcal{I}_n & & & W_s
        \end{bmatrix}\begin{bmatrix}
            \svec(\Delta S)\\
            \svec(\Delta Z)\\
            \svec(\Delta H)\\
            \svec(\Delta Q)
        \end{bmatrix}=\begin{bmatrix}
            0\\
            0\\
            0\\
            b
        \end{bmatrix},
\end{equation}
where $b=\svec(\mu Q^{-1}-S)$, $E=U_{n}([I_n~0]\otimes [I_n~0])U_{n+m}^{\top}$, $G=\hat{\Phi}_t \otimes_s\hat{\Phi}_t-\mathcal{I}_{n+m}$, $W_s=W^{-1}\otimes_sW^{-1}$, $U_{n}$ and $U_{n+m}$ are defined in Appendix \ref{appendix:svex},  and
    \begin{align*}
        &F=(\begin{bmatrix}
        0_{n\times m
        } \\ I_m
    \end{bmatrix}\begin{bmatrix}
        0_{m\times n} & I_m
    \end{bmatrix})\otimes_s(\begin{bmatrix}
        I_n \\ -\bar{K}
    \end{bmatrix}\begin{bmatrix}
        I_n & -\bar{K}^{\top}
    \end{bmatrix}).
    \end{align*}
\end{proposition}

Now we are ready to define the primal-dual interior-point iterations. Starting at a strictly feasible point with the search direction computed by \eqref{eq:newton_svec}, we consider the full NT-step
\begin{equation}\label{eq:NTupdates}
\begin{aligned}
    & Q_+= Q + \Delta Q, \quad S_+= S + \Delta S,   \\
    & Z_+= Z + \Delta Z, \quad  H_+= H + \Delta H, 
\end{aligned}
\end{equation}
to eliminate the need to compute damped step sizes as in the line search framework.  Reduce $\mu$ by $\mu_+<\mu$  and continue the iterations. As will be discussed in the next part, the primal-dual interior point method will finally converge as $\mu \rightarrow 0$.

\subsection{Convergence Analysis of the Offline Algorithm}\label{sec:SDP_convergence}
In this part, convergence of the full NT-step primal-dual interior-point iterations in \eqref{eq:NTupdates} is studied. To facilitate analysis, we define the scaled search direction:
\begin{equation}\label{eq:scale_D}
\begin{aligned}
&D_Q= \tfrac{1}{\sqrt{\mu}} W^{-\tfrac{1}{2}}\Delta QW^{-\tfrac{1}{2}}, \, D_S= \tfrac{1}{\sqrt{\mu}}W^{\tfrac{1}{2}}\Delta SW^{\tfrac{1}{2}}, \\
&D_Z=  \tfrac{1}{\sqrt{\mu}}\tilde{W}\Delta Z\tilde{W}, \, D_H=  \tfrac{1}{\sqrt{\mu}} \tilde{W}^{-1}\Delta H\tilde{W}^{-1},
\end{aligned}
\end{equation}
where $\tilde{W}=\diag (W^{\tfrac{1}{2}},I_m)$. Therefore, the search direction defined in (\ref{eq:newton}a-\ref{eq:newton}c) and \eqref{eq:NT_QS} is equivalent to
\begin{subequations}\label{eq:modified_newton}
    \begin{align}
        &D_S=\begin{bmatrix}
        I_n & 0
    \end{bmatrix}D_Z\begin{bmatrix}
        I_n \\ 0
    \end{bmatrix} \label{eq:DS},\\
   &D_Z=\overline{\tilde{\Psi}(D_H)} + \tilde{\Phi}_{t} D_Z \tilde{\Phi}^{\top}_{t}, \label{eq:DZ}\\
   &D_H=\diag(D_Q, 0)+\tilde{\Phi}_{t}^{\top}D_H{\tilde{\Phi}_{t}}, \label{eq:DH}\\
   &D_Q+D_S=V^{-1}-V, \label{eq:DQ}
    \end{align}
\end{subequations}
where $\tilde{\Phi}_{t}=\tilde{W}\hat{\Phi}_t\tilde{W}^{-1}$ and
\begin{align*}
\tilde{\Psi}(D_{H})=\begin{bmatrix}0 \\ I_m\end{bmatrix}\begin{bmatrix}0 & I_m\end{bmatrix}D_H\begin{bmatrix}W^{\tfrac{1}{2}} \\ -\bar{K}\end{bmatrix}\begin{bmatrix}W^{\tfrac{1}{2}} & -\bar{K}^{\top}\end{bmatrix}.
\end{align*}

Define the centrality function as
\begin{equation}\label{eq:def_delta}
\delta(V):=\delta(Q,S;\mu)=\norm{I-V^2}.
\end{equation}
Note that $QS=\mu I_n$ if and only if $\delta(V)=0$. In essence, $\delta(V)$ characterizes the current deviation from the central path.

\begin{remark}
Different from most of the existing literature that uses $\delta(Q,S;\mu)=\norm{D_Q+D_S}$, namely, $\delta(V)=\norm{V^{-1}-V}$ as in our case, we adopt \eqref{eq:def_delta} to facilitate the convergence analysis of the online algorithm in Section \ref{sec:online}. To this end, convergence of the full NT-step algorithm with offline data using such a centrality function is first studied as follows, which will serve as the foundation for the online iterations in the next part.
\end{remark}

The main challenge of generalizing the quadratic convergence of $\delta(V)$ in linear SDP to nonlinear cases lies in that the orthogonality of scaled search directions no longer holds, which is studied for our problem in Lemmas \ref{lem:Tr_DsDq} and \ref{lem:spectral_DqDs}. In the following Lemmas, we define $P_V=D_Q+D_S$, $Q_V=D_Q-D_S$, $D_{QS}= \frac{1}{2}(D_Q D_S+D_SD_Q)$. 
\begin{lemma}\label{lem:Tr_DsDq}
Let $(D_{S},D_{Z},D_{H},D_{Q})$ solve \eqref{eq:modified_newton} and assume $\delta<1$, then it holds that
    \begin{equation}\label{eq:Tr_DQDS}
        0\leq \langle D_S, D_Q \rangle \leq \frac{\delta^2}{4(1-\delta)}.
    \end{equation}
\end{lemma}
\proof
Since $D_Z$ and $D_H$ can be explicitly solved from the Lyapunov equations \eqref{eq:DZ} and \eqref{eq:DH}, we have
\begin{align*}
         \langle D_S, D_Q \rangle &= \tr (D_Z \begin{bmatrix}
             I_n \\ 0
         \end{bmatrix}D_Q\begin{bmatrix}
             I_n & 0
         \end{bmatrix}) \\
         &=\tr(\sum_{i=0}^{\infty}{\tilde\Phi_{t}^i}\overline{\tilde\Psi(D_H)}{(\tilde\Phi_{t}^i)}^{\top} \begin{bmatrix}
             D_Q & 0\\ 0 & 0
         \end{bmatrix})\\      &=\tr(\overline{\tilde\Psi(D_H)}\sum_{i=0}^{\infty}{(\tilde\Phi_{t}^i)}^{\top} \begin{bmatrix}
             D_Q & 0\\ 0 & 0
         \end{bmatrix}{\tilde\Phi_{t}^i})\\
         &=\tr(\overline{\tilde\Psi(D_H)}D_H)\\
         &=\norm{\begin{bmatrix}
        0_{m\times n} & I_m
    \end{bmatrix}D_H\begin{bmatrix}
        W^{1/2} \\ -\bar{K}
    \end{bmatrix}}^2\geq 0.
\end{align*}    
In addition, we have $V^2=I_{n}-(I_{n}-V^2)\succeq (1-\delta(V))I_{n}$, which implies
\begin{equation}\label{eq:norm_DQ+DS}
\begin{aligned}
\norm{P_{V}}^{2} = & \norm{V^{-1}-V}^2  \leq \norm{V^{-1}}^2_2\norm{I_{n}-V^2}^2\leq \frac{\delta^2}{1-\delta}.
\end{aligned}        
\end{equation} 
    On the other hand, it holds that 
\begin{align*}
\norm{P_{V}}^{2} = \norm{Q_{V}}^{2} + 4\langle D_S, D_Q \rangle \geq 4\langle D_S, D_Q \rangle,
\end{align*}   
which implies the right-hand side of \eqref{eq:Tr_DQDS}.
\endproof


\begin{lemma}\label{lem:spectral_DqDs}
Let $D_{Q}, D_{S} \in \mathbb{S}^n_+$ satisfy $\langle D_Q, D_S \rangle \geq 0$. The spectral radius of $D_{QS}$ is bounded via $\delta$ by
\begin{equation}
\rho(D_{QS}) \leq \frac{\delta^{2}}{4(1-\delta)}.
\end{equation}
\end{lemma}
\proof
Since $\langle D_Q, D_S \rangle \geq 0$, we have $\norm{Q_{V}}^{2} \leq \norm{P_{V}}^{2}$. 
By rewriting $D_{QS} = \frac{1}{4}(P_{V}^{2} - Q_{V}^{2})$, we further have
\begin{equation*}
\begin{aligned}
\frac{1}{4}\norm{P_{V}}^{2} I_{n} \succeq D_{QS} \succeq
-\frac{1}{4} \norm{Q_{V}}^{2} I_{n} \succeq -\frac{1}{4} \norm{P_{V}}^{2} I_{n},
\end{aligned}
\end{equation*}
which yields $\rho(D_{QS}) \leq \frac{1}{4} \norm{P_{V}}^{2} \leq \frac{\delta^{2}}{4(1-\delta)}$.
\endproof

Next, we show that the full NT-step guarantees a strict feasible update in Proposition \ref{prop:feasible_step} and quadratic convergence to the central path in Proposition \ref{prop:quad_converge}.
\begin{proposition}[Strict feasibility of the full NT-step]\label{prop:feasible_step}
    Assume $Q,S\succ 0$ and $\delta<2\sqrt{2}-2$, then $Q_+,S_+\succ 0$.
\end{proposition}
\proof
The proof follows a similar vein as that for the standard linear SDP \cite[Lemma 7.1]{Klerk2002}, which in essence requires $\rho(D_{QS}) < 1$. Hence, the claim can be proved by Lemma \ref{lem:spectral_DqDs}. 
\endproof

\begin{proposition}[Local quadratic convergence] \label{prop:quad_converge}
    Assume $\delta<1$, then 
    \begin{equation}
        \delta_+:=\delta(Q_+,S_+;\mu)\leq\frac{\delta^2}{2\sqrt{2}(1-\delta)}.
    \end{equation}
    Specifically, if $\delta\leq\frac{1}{2}$, then $\delta_+\leq\frac{1}{\sqrt{2}}\delta^2$.
\end{proposition}    
\proof
By \eqref{eq:norm_DQ+DS}, we have
\begin{align*}
\norm{D_{QS}}^2&=\norm{\frac{1}{4}(P_V^2-Q_V^2)}^2\leq\frac{1}{16}(\norm{P_V}^4+\norm{Q_V}^4)\\
&\leq\frac{1}{8}\norm{P_V}^4\leq\frac{\delta^4}{8(1-\delta)^2}.
\end{align*}
Let $V_+$ define the spectral geometric mean of $S_+$ and $Q_+$ as in \eqref{eq:def_V}. Similar to Lemma 7.4 in \cite{Klerk2002}, one can show that
\begin{align*}
        \delta_+^2=\norm{I_{n}-V_+^2}^2\leq\norm{D_{QS}}^2\leq\frac{\delta^4}{8(1-\delta)^2}.
\end{align*}
If $\delta\leq\frac{1}{2}$, it is clear that $        \delta_+\leq\frac{\delta^2}{2\sqrt{2}(1-\delta)}\leq \frac{\delta^2}{\sqrt{2}}$. 
\endproof 

Finally, convergence of the full-NT steps is studied by incorporating the influence of reducing $\mu$. Through alternating between \eqref{eq:NTupdates} and reducing $\mu$ by $\mu_+=\theta\mu$ with $0<\theta<1$, we show in Theorem \ref{thm:delta_neighbor_offline} that the iterates will always remain strictly feasible and lie in the neighborhood of the central path by $\delta\leq\frac{1}{2}$, thus leading to a vanishing duality gap as $\mu\rightarrow0$. It is noteworthy that in Theorem \ref{thm:delta_neighbor_offline} we give a tighter bound of $\delta(Q_+,S_+;\mu_+)$, namely, $c\leq\frac{1}{2}$, which will help the convergence analysis of online iterations in the next part. 
\begin{theorem}\label{thm:delta_neighbor_offline}
    For any $c\in(0,\frac{1}{2}]$, let $\mu_+=\theta\mu$ with $\theta$ satisfying $\mathcal{A}(c)\leq\theta<1$, where
    \begin{equation}\label{eq:theta_bnd}
        \mathcal{A}(c):=\frac{8n+1-\sqrt{(64c-2)n+18c+1}}{8(n-c)} .
    \end{equation}
    If $\delta\leq\frac{1}{2}$, then the full NT-step is strictly feasible and the updated centrality function is bounded by $\delta(Q_+,S_+;\mu_+)\leq c$.
\end{theorem}
\proof
Strict feasibility follows by Proposition \ref{prop:feasible_step}. By \eqref{eq:def_V}, we know that $V_+\gets \frac{1}{\sqrt{\theta}}V_+$ if $\mu_+\gets \theta\mu$. Then we have
\begin{align}
&\delta(Q_+,S_+;\mu_+)=\norm{I_{n}-\frac{V_+^2}{\theta}}^2=\norm{\frac{I_{n}-V_+^2}{\theta}-\frac{(1-\theta) I_{n}}{\theta}}^2 \notag \\=&\frac{1}{\theta^2}\delta_+^2+\frac{n(1-\theta)^2}{\theta^2}-\frac{2(1-\theta)}{\theta^2}\langle I_{n}-V_+^2,I_{n} \rangle \label{eq:delta_mu+} \\
=&\frac{1}{\theta^2}\delta_+^2+\frac{n(1-\theta)^2}{\theta^2}+\frac{2(1-\theta)}{\theta^2}(\norm{V_+}^2-n). \notag
\end{align}
By Proposition \ref{prop:quad_converge}, $\delta\leq\frac{1}{2}$ implies $\delta_+\leq\frac{1}{\sqrt{2}}\delta^2\leq \frac{1}{4\sqrt{2}}$. Using \eqref{eq:def_V} and \eqref{eq:scale_D}, it holds that
    \begin{align*}
        &Q_+=Q+\Delta Q=\sqrt{\mu}W^{\frac{1}{2}}(V+D_Q)W^{\frac{1}{2}},\\
        &S_+=S+\Delta S=\sqrt{\mu}W^{-\frac{1}{2}}(V+D_S)W^{-\frac{1}{2}}.
    \end{align*}
    Then the duality gap is bounded via Lemma \ref{lem:Tr_DsDq} by
\begin{equation}\label{eq:duality_gap_Q+S+}
\begin{aligned}
& \langle Q_+,S_+ \rangle =\mu\tr[(V+D_Q)(V+D_S)]\\
= & \mu[\langle V^{-1},V\rangle+\langle D_Q,D_S\rangle] \leq\mu(n+\frac{\delta^2}{4(1-\delta)}),
\end{aligned}
\end{equation}
where the second equality follows the fact that $D_Q+D_S=V^{-1}-V$.
Since $\frac{\delta^2}{1-\delta}$ monotonically increases with $\delta\in[0,\frac{1}{2}]$, \eqref{eq:def_V} further yields
    \begin{equation*}
        \norm{V_+}^2=\frac{1}{\mu}\langle Q_+,S_+ \rangle  \leq n+\frac{\delta^2}{4(1-\delta)}\leq n+\frac{1}{8}.
    \end{equation*}
    Substituting the above results in \eqref{eq:delta_mu+} gives
    \begin{align*}
        \delta(Q_+,S_+;\mu_+)\leq\frac{1}{\theta^2}(\frac{1}{32}+n(1-\theta)^2+\frac{1-\theta}{4}),
    \end{align*}
    where the right-hand side monotonically decreases with $\theta\in(0,1)$. Then straightforward computations yield that \eqref{eq:theta_bnd} implies $\delta(Q_+,S_+;\mu_+)\leq c$.
\endproof

\subsection{Adaptive Data-Driven IRL with Online Data}\label{sec:online}  
In this part, the online adaptation recursions are developed to solve a time-varying variant of \eqref{eq:Jt}, where the primal-dual iterates are coupled with the online collection of system data $(x_t,u_t)$. At each time $t$, a modified one-step primal-dual interior point iterate is applied based on $(X_{0,t},U_{0,t},X_{1,t})$, following which new system data  $(x_{t+1},u_{t+1})$ is collected at time $t+1$ and the new search direction is computed.

\subsubsection{Outline of the Online Adaptive Algorithm}
Compared with the offline data-driven framework in \eqref{eq:Jt}, the main challenge of the online iterative algorithm lies in that the cost function $f_t$ is time-varying as the iterations continue. 
More specifically, the time-varying property of $\hat\Phi_t$ implies that the Newton's step applied to \eqref{subeq:Z_eqn} and \eqref{subeq:H_eqn} is no longer consistent with \eqref{eq:newton_z} and \eqref{eq:newton_H}. To this end, given $(Q_t,S_t,Z_t,H_t)$ that satisfies the IPC at iteration $t$, a new feasible search direction is derived via the following two steps such that $(Q_{t+1},S_{t+1},Z_{t+1},H_{t+1})$  satisfies the IPC at iteration $t+1$. Firstly, since $\hat{\Phi}_{t+1}$ is stable under the assumption in \eqref{eq:SNR_bnd}, we define $(\bar{S}_t,\bar{Z}_t,\bar{H}_t)$ as the unique solution to the following equations at $Q_t \succ 0$:
\begin{equation}\label{eq:recur_t_inter}
\begin{aligned}
& \bar{S}_t=\begin{bmatrix}I_n & 0\end{bmatrix}\bar{Z}_t\begin{bmatrix}I_n \\ 0\end{bmatrix},\\
&\bar{Z}_t=\overline{\Psi(\bar{H}_{t})} + \hat{\Phi}_{t+1} \bar{Z}_t\hat{\Phi}_{t+1}^{\top},\\
&\bar{H}_t=\diag(Q_{t}, I_{m})+\hat{\Phi}_{t+1}^{\top}\bar{H}_t\hat{\Phi}_{t+1}. \end{aligned}
\end{equation}
If $\bar{S}_{t}\succ0$, there is a unique weighting point $\bar{W}_t\succ 0$ such that $\bar{W}_{t}^{-1}Q_t\bar{W}_{t}^{-1}=\bar{S}_{t}$. Similar to Proposition \ref{prop:NT_direction}, there exists unique $(\Delta S_{t},\Delta Z_{t},\Delta H_{t},\Delta Q_{t})$ that solves
\begin{equation}\label{eq:DQ_online}
\begin{aligned}
&\Delta S_t=\begin{bmatrix}I_n & 0\end{bmatrix}\Delta Z_t\begin{bmatrix}I_n \\ 0\end{bmatrix}, \\
&\Delta Z_t=\overline{\Psi(\Delta H_t)} + \hat{\Phi}_{t+1} \Delta Z_t\hat{\Phi}_{t+1}^{\top}, \\
&\Delta H_t=\diag(\Delta Q_t, 0)+\hat{\Phi}_{t+1}^{\top}\Delta H_t\hat{\Phi}_{t+1}, \\
&\bar{W}_t^{-1}{\Delta Q_t}\bar{W}_t^{-1}+\Delta S_t =\mu Q_t^{-1}-\bar{S}_t.
    \end{aligned}
\end{equation}
Then the primal and dual variables are updated by
\begin{equation}\label{eq:Q+_online}
\begin{aligned}
& Q_{t+1}=Q_t+\Delta Q_t, \quad S_{t+1}=\bar{S}_t+\Delta S_t, \\
& Z_{t+1}=\bar{Z}_t+\Delta Z_t, \quad H_{t+1}=\bar{H}_t+\Delta H_t.
\end{aligned}
\end{equation}
which obviously satisfy the IPC at $t+1$.

If we can show that $\bar\delta_t(\bar V):=\delta(Q_t,\bar{S}_t;\mu_t)\leq\frac{1}{2}$, then following the convergence results of the offline algorithm in Section \ref{sec:SDP_convergence}, the full NT-step \eqref{eq:DQ_online} will ensure a vanishing dualtiy gap. The online algorithm is outlined in Algorithm \ref{alg:Newton_data} and its convergence analysis is given in the next part.
\begin{algorithm}[t]
\caption{Adaptive online data-driven IRL.}\label{alg:Newton_data}
\begin{algorithmic}
    \REQUIRE 
    \STATE initial dataset $(X_{0,t_0},U_{0,t_0},X_{1,t_0})$;
    \STATE strictly feasible $(Q_{t_0},H_{t_0},Z_{t_0},S_{t_0})$ and $\mu_{0} = \frac{\langle Q_{t_0},S_{t_0}\rangle }{n}$ such that $\delta(Q_{t_0},S_{t_0};\mu_{0}) \leq \frac{1}{2}-h_7 $ .
    \ENSURE {adaptive estimate of the cost matrix $Q_t$. }
    \WHILE{$t>t_0$} 
        \STATE Let $\mu={\mu_0}/{(t-t_0)}$.
        \STATE Collect new data $(x_t, u_t)$.
        \STATE Update $(\bar{X}_{1,t},\Lambda_{t}^{-1})$ recursively by \eqref{eq:data_recur} to obtain $\hat{\Phi}_t$.
        \STATE Compute $(\bar{S}_{t-1},\bar{Z}_{t-1},\bar{H}_{t-1})$ by solving \eqref{eq:recur_t_inter}.
        \STATE Solve $(\Delta S_{t-1}, \Delta Z_{t-1},\Delta H_{t-1},\Delta Q_{t-1})$ by \eqref{eq:DQ_online} and update $(S_{t},Z_{t},H_{t},Q_{t})$ with \eqref{eq:Q+_online}.
    \ENDWHILE
\end{algorithmic}
\end{algorithm}

\subsubsection{Convergence Analysis}
In the online optimization framework, at each time a new central path is defined and local convergence is required to be guaranteed in the presence of process noise.
For brevity, in the sequel we denote $\Sigma_{\Phi_{\bar{K}}} = \sum_{i = 0}^{\infty} \Phi_{\bar{K}}^{i}(\Phi_{\bar{K}}^{i})^{\top}$, $\Sigma_{\hat{\Phi}_{t}} = \sum_{i = 0}^{\infty} \hat{\Phi}_{t}^{i}(\hat{\Phi}_{t}^{i})^{\top}$ and
\begin{align*}
&h_{1} = \norm{\bar{K}}_2+1, \, h_{2} = 1+\norm{\Phi_{\bar{K}}}_{2}, \\ & h_{3} = \max\lbrace \norm{\Sigma_{\Phi_{\bar{K}}}}_{2}, \norm{\Sigma_{\Phi^{\top}_{\bar{K}}}}_{2} \rbrace.
\end{align*}

Since the time-varying property mainly results from the updates of $\hat{\Phi}_t$, its variation is first analyzed in Proposition \ref{prop:Phi_disturb}, whose proof is given in Appendix \ref{append:Phi_disturb}.
\begin{proposition}\label{prop:Phi_disturb}
    Assume $\lbrace x_t\rbrace_{t\geq 0}$ and $\lbrace u_t\rbrace_{t\geq 0}$ are upper bounded. Then the variation of $\hat{\Phi}_t$ is bounded by
    \begin{equation}
        \begin{aligned}
            &\norm{\hat{\Phi}_{t}-\Phi_{\bar{K}}}_{2} \leq \frac{2h_{1}\kappa}{\gamma\zeta},\\
            &\norm{\Delta_{\hat{\Phi},t}}_{2} :=\norm{\hat{\Phi}_{t+1} - \hat{\Phi}_{t}}_{2}\leq \frac{h_4}{t},
        \end{aligned}
    \end{equation}
    where $h_4>0$ is some constant that only depends on the system characteristics.
\end{proposition}

Next, the Slater regularity condition is strengthened to take into account the variation of $\hat{\Phi}_t$. Later we will see that Assumption \ref{assump:perturbed_slater} in essence guarantees the strict feasibility of $f_t$ for any $t\geq t_0$.
\begin{assumption}[Perturbed Slater regularity condition]\label{assump:perturbed_slater}
    Suppose there exists $(\tilde{Q},\tilde{S},\tilde{H},\tilde{Z})$ that satisfies the IPC at $t_0$ and $\tilde{S}\succ \tau_s (1 + \norm{\tilde{Q}}_{2})I_n$, where
\begin{equation}
    \tau_s= 128 h_{1}^{3}h_{3}^{3}(h_{2} + h_{1}\frac{2\kappa}{\gamma \zeta}) \frac{2\kappa}{\gamma \zeta}. 
\end{equation}
\end{assumption}

For brevity, in the sequel we define
\begin{equation*}
    h_5=\langle \tilde{Q},\tilde{S} \rangle +\tau_s (1 + \norm{\tilde{Q}}_{2})\tr(\tilde{Q}) -\frac{\mu_{0}}{t_0}\log\det(\tilde{Q}),
\end{equation*}
if $\det(\tilde{Q})\leq 1$, while 
$h_5=\langle \tilde{Q},\tilde{S} \rangle +\tau_s (1 + \norm{\tilde{Q}}_{2})\tr(\tilde{Q})$ if $\det(\tilde{Q})> 1$, and
\begin{align*}
 &h_6=-\frac{\mu_0}{\tilde{\sigma}_st_0}W_{-1}(-\frac{\tilde{\sigma}_st_0}{\mu_0}\exp(-{\frac{h_4t_0}{\mu_0}})), \\
    &h_7=\frac{64\sqrt{n}}{\mu_{0}}h_{1}^{2}h_{3}^{3}h_4(h_{2} + h_{1}\frac{2\kappa}{\gamma \zeta})  (1 + h_6)h_6,
\end{align*}
where  $W_{-1}(\cdot)$ denotes the Lambert $W$-function and $\tilde{\sigma}_s$ is defined as $\tilde{\sigma}_s=\underline{\sigma}(\tilde{S})-\tau_s (1 + \norm{\tilde{Q}}_{2})$.

To further analyze the online updates of the central path, we give Lemma \ref{lemma:lyapunov_disturbation} on the perturbation of Lyapunov equations,  whose proof follows a similar vein of \cite{Fazel2018Global} and is omitted due to page limitations.
\begin{lemma}\label{lemma:lyapunov_disturbation}
Consider the discrete-time Lyapunov equation $\Sigma=Q+A\Sigma A^{\top}$ where $A$ is stable. Let $A$ be perturbed to some stable matrix $\bar{A}=A+\Delta$ and denote $\bar \Sigma$ as the corresponding solution. Define $\Sigma_I:=\sum_{i=0}^{\infty}A^i(A^i)^{\top}$. If $\norm{\Delta}_2 \leq 1/(4\norm{\Sigma_I}_2(1+\norm{A}_2))$, then $\norm{\bar{\Sigma}-\Sigma}_2\leq 4\norm{Q}_2\norm{\Sigma_I}_2^2(1+\norm{A}_2)\norm{\Delta}_2$.
\end{lemma}

With the above results, we next show that the online iterations are always bounded within a neighborhood of the time-varying central path as long as the SNR is large enough, which consequently guarantees the convergence of the modified full-NT steps. 
\begin{theorem}\label{thm:convergence_online}
Given $\mu_0>0$ and let $\mu_t=\frac{\mu_0}{t-t_0}$. Assume $\frac{\gamma}{\kappa} \geq \max \{\frac{2\omega}{\zeta} , \frac{2\hat{\sigma}}{\zeta}, \tau_{2}, \tau_{3} \}$, where
\begin{equation}
\begin{aligned}
\tau_{1} &= \tfrac{8}{\zeta}h_{1}h_{2}h_{3},\\
\tau_{2} & = \frac{16h_1h_3h_4}{\zeta(t_0-8h_2h_3h_4)}, \\
\tau_{3} & = 2\tau_{1} + \frac{8}{\zeta}h_{1}\sqrt{h_{3} + 4h_{2}^{2}h_{3}^{2}}.
\end{aligned}
\end{equation}
Suppose $\mu_0$ and $\frac{\gamma}{\kappa}$ are chosen such that $h_7<\frac{1}{2}$ and $\mathcal{A}(\frac{1}{2}-h_7)\leq\frac{1}{2}$.
Then if $\delta_t\leq\frac{1}{2}-h_7$ , it holds that $\delta_{t+1} \leq \frac{1}{2}-h_7$.
\end{theorem}
\begin{remark}
    The threshold for SNR only depends on the intrinsic properties of the nominal system. In case of lacking such prior knowledge, a conservative approximation could be obtained from the initial dataset (such as the polytopic case in \cite{LiPolytopic2024}) at the cost of requiring a larger SNR.
\end{remark}
\noindent\hspace{2em}{\itshape Proof of Theorem \ref{thm:convergence_online}: }
The rationale for analyzing the online algorithm is to incorporate the perturbation analysis of two cascaded Lyapunov functions and the convergence results in Section \ref{sec:SDP_convergence}. Below we outline the main proof steps, where the proofs of the lemmas are given in the Appendix.

\noindent\textbf{(i) Upper bounds for $\norm{\Sigma_{\hat{\Phi}_{t}}}_{2}$.} Since $\frac{\gamma}{\kappa} \geq \tau_3> \tau_{1}$, we can obtain that
\begin{equation*}
\begin{aligned}
\norm{\Phi_{\bar{K}} - \hat{\Phi}_{t}}_{2}  \leq h_{1} \frac{2\kappa}{\gamma \zeta} \leq \frac{2h_{1        }}{\zeta} \frac{1}{\tau_{1}}  \leq \frac{1}{4h_{2}h_{3}}.
\end{aligned}
\end{equation*}
Then by Lemma \ref{lemma:lyapunov_disturbation}, we have 
\begin{equation}\label{eq:Sigma_Phit_bnd}
    \norm{\Sigma_{\hat{\Phi}_{t}}}_{2}  \leq \norm{\Sigma_{\Phi_{\bar{K}}}}_{2} + \norm{\Sigma_{\Phi_{\bar{K}}} - \Sigma_{\hat{\Phi}_{t}}}_{2}  \leq 2h_{3}.
\end{equation}
\textbf{(ii) Upper bounds for $\norm{H_{t} - \bar{H}_{t}}_{2}$ and $\norm{H_t}_2$.}
\begin{lemma}\label{prop:H_perturb_bnd}
    If $\frac{\gamma}{\kappa} \geq\tau_{2}$, it holds that
\begin{equation}
\norm{H_{t} - \bar{H}_{t}}_{2} \leq \frac{16}{t}h_{3}^{2}h_{4}(h_{2} + h_{1}\frac{2\kappa}{\gamma \zeta})  (1 + \norm{Q_{t}}_{2}).
\end{equation}

Furthermore, $\norm{H_{t}}_{2} \leq 2(1+\norm{Q_{t}}_{2})h_{3}$. 
\end{lemma}
\textbf{(iii) Upper bound for $\norm{Z_{t} - \bar{Z}_{t}}_{2}$. }
\begin{lemma} \label{proposition:normbound_Zt}
If $\frac{\gamma}{\kappa} \geq \tau_{2}$, it holds that
\begin{equation*}
\norm{Z_{t} - \bar{Z}_{t}}_{2} \leq \frac{64}{t} h_{1}^{2}h_{3}^{3}h_4(h_{2} + h_{1}\frac{2\kappa}{\gamma \zeta})  (1 + \norm{Q_{t}}_{2}).
\end{equation*}
\end{lemma}
\textbf{(iv) Uniform bound for $\norm{Q_{t}}_2$. }  
We first show that under Assumption \ref{assump:perturbed_slater}, the Slater regularity condition is satisfied at any $t\geq t_0$. Denote $(\hat{S}_t,\hat{H}_t,\hat{Z}_t)$ as the unique matrix triple that satisfies \eqref{subeq:S_eqn}, \eqref{subeq:Z_eqn}, \eqref{subeq:H_eqn} with given $\tilde{Q}$. It is clear that $\norm{\hat{\Phi}_{t}-\hat{\Phi}_{t_0}}_{2}\leq\norm{\hat{\Phi}_{t_0}-\Phi_{\bar{K}}}_{2}+\norm{\hat{\Phi}_{t}-\Phi_{\bar{K}}}_{2}\leq\frac{4h_{1}\kappa}{\gamma\zeta}$. Under the assumption that $\frac{\gamma}{\kappa} \geq \tau_3$, following similar proof as in Lemmas \ref{prop:H_perturb_bnd} and \ref{proposition:normbound_Zt}, one can show that
\begin{equation}\label{eq:hatS_bnd}
\begin{aligned}
    &\norm{\hat{S}_{t} - \tilde{S}}_{2} \leq \norm{\hat{Z}_{t} - \tilde{Z}}_{2}
    \leq \tau_s (1 + \norm{\tilde{Q}}_{2}).
\end{aligned}
\end{equation}
Hence, $\tilde{S}\succ \tau_s (1 + \norm{\tilde{Q}}_{2})I_n$ implies that $\hat{S}_{t}\succ 0$, namely, $(\tilde{Q}, \hat{S}_t,\hat{H}_t,\hat{Z}_t)$ is strictly feasible for $f_{\mu}(Q)$ at any $t>t_0$. Therefore, it is natural to consider the iterations in the level set  $\lbrace Q\succ 0: f_{\mu}(Q) \leq f_{\mu}(\tilde{Q}) \rbrace$, otherwise $\tilde{Q}$ provides a better solution. Denote $\lambda$ as an arbitrary eigenvalue of $Q_t$. Similar to \eqref{eq:eigQ_bound}, by \eqref{eq:hatS_bnd} we have
\begin{equation*}
    \begin{aligned}
       \tilde{\sigma}_s\lambda-\mu\log\lambda\leq &f_{\mu}(\tilde{Q})+\frac{1}{2}	\langle \hat{H}_t,\Psi(\hat{H}_t) \rangle-\tr(\diag(0, I_{m})\hat{Z}_t)\\
        =& \langle \tilde{Q},\hat{S}_t \rangle-\mu\log\det(\tilde{Q})\leq h_5,
    \end{aligned}
\end{equation*}
which implies
    $\lambda \leq -\frac{\mu}{\tilde{\sigma}_s}W_{-1}(-\frac{\tilde{\sigma}_s}{\mu}\exp(-{\frac{h_5}{\mu}})).$
Since the right-hand side monotonically increases with $\mu$, we then have $\norm{Q_t}_2\leq h_6$.


\noindent\textbf{(v) Upper bound for $\bar{\delta}_{t}$. }  Finally, we can show that
\begin{equation}
\begin{aligned}
\bar{\delta}_{t} & = \norm{I_{n} - \bar{V}_{t}^{2}} \leq \norm{I_{n} - V_{t}^{2}} + \norm{V_{t}^{2} - \bar{V}_{t}^{2}} \\ &  = \delta_{t} + \frac{1}{\mu_{t}} \norm{Q_{t}^{\frac{1}{2}}(S_{t} - \bar{S}_{t}) Q_{t}^{\frac{1}{2}}   } \\ & \leq \delta_{t} + \frac{\sqrt{n}}{\mu_{t}} \norm{Q_{t}}_{2} \norm{S_{t} - \bar{S}_{t}}_{2}.
\end{aligned}
\end{equation}
Since $\norm{S_{t} - \bar{S}_{t}}_{2} \leq \norm{Z_{t} - \bar{Z}_{t}}_{2}$,  Lemma \ref{proposition:normbound_Zt} yields
\begin{equation*}
\begin{split}
\bar{\delta}_{t} & \leq \delta_{t}+\frac{64\sqrt{n}}{t\mu_{t}}h_{1}^{2}h_{3}^{3}h_4(h_{2} + h_{1}\frac{2\kappa}{\gamma \zeta})  (1 + \norm{Q_{t}}_{2})\norm{Q_{t}}_{2} \\ 
& \leq \delta_{t}+h_7\leq\frac{1}{2}.
\end{split}
\end{equation*}
Furthermore, we have $\frac{\mu_{t+1}}{\mu_{t}}=\frac{t-t_0}{t-t_0+1}\geq\frac{1}{2}\geq\mathcal{A}(\frac{1}{2}-h_7)$. By Theorem \ref{thm:delta_neighbor_offline} we can directly obtain that $\delta_{t+1} \leq \frac{1}{2}-h_7$.
\endproof

In conclusion, we have shown that if $\delta_{t_0}\leq \frac{1}{2}-h_7$, the online iterations will always stay in the same neighborhood near the time-varying central path, namely, $\delta_t\leq \frac{1}{2}-h_7$ for any $t\geq t_0$. The duality gap will then converge to zero in a sublinear rate as shown in Proposition \ref{prop:complexity}.

\begin{proposition}[Non-asymptotic convergence]\label{prop:complexity}
Choose $(Q_{t_0},S_{t_0})$ be strictly feasible and $\mu_{0} = {\langle Q_{t_0},S_{t_0}\rangle}/{n}$ such that $\delta(Q_{t_0},S_{t_0};\mu_{0}) \leq \frac{1}{2}-h_7$. Under the assumptions in Theorem \ref{thm:convergence_online}, let $Q_{k}$ and $S_{k}$ denote the solutions after $k$ iterations. Then the duality gap converges in a sublinear rate with order $\mathcal{O}(1/k)$, and $\langle Q_{k},S_{k} \rangle \leq \epsilon$ holds for
\begin{equation}\label{eq:iteration_bounds}
k \geq \frac{9}{8\epsilon}\langle Q_{t_0},S_{t_0} \rangle .
\end{equation}
\end{proposition}
\proof
In Theorem \ref{thm:convergence_online} we have shown that $\bar{\delta}_{t}\leq\frac{1}{2}$ for any $t\geq t_0$. By \eqref{eq:duality_gap_Q+S+}, it follows that
\begin{equation*}
\langle Q_{t},S_{t} \rangle \leq \mu_t (n + \frac{\bar{\delta}_{t-1}^{2}}{4(1-\bar{\delta}_{t-1})}) \leq \mu_t(n + \frac{1}{8}) \leq \frac{9}{8}n\mu_t.
\end{equation*}
Denote $k:=t-t_0$ as the number of iterations. Then we can obtain that
\begin{equation*}
\langle Q_{t},S_{t}\rangle \leq \frac{9}{8}n\mu_{t} = \frac{9n\mu_{t_0}}{8k} = \frac{9}{8k}\langle Q_{t_0},S_{t_0}\rangle.
\end{equation*}
Thus, $\langle Q_{t},S_{t}\rangle \leq \epsilon$ is achieved if $\frac{9}{8k}\langle Q_{t_0},S_{t_0}\rangle \leq \epsilon$, which implies $k \geq \frac{9}{8\epsilon}\langle Q_{t_0},S_{t_0} \rangle$.
\endproof

\begin{remark}
Recall that on the central path $QS=\mu I_n$ and the solution space of $Q$ to the IRL problem is bounded in $\mathbb{S}^n_{++}$. Then $\hat{S}_t$ and $\Delta_{\Phi,t}$ converge in the same order by $\mathcal{O}(1/k)$. Hence, given $S_t\succ 0$, it is easy to satisfy $\bar{S}_t\succ 0$ with a large but bounded SNR, which, otherwise, can be satisfied by choosing a step size $\alpha\in(0,1)$ and all the properties in Section \ref{sec:SDP_convergence} still hold. 
\end{remark}

\begin{remark}
    In model-free implementation of the IRL algorithm, the proposed adaptive learning rule does not require the expert (or a learner with the same dynamics) to collect online system data in a pre-defined manner as existing IRL results. We rigorously show how arbitrary PE data can help improve the cost learning accuracy only under mild assumptions on the SNR. Non-asymptotic convergence is studied in the presence of system uncertainties, which has been largely unexploited in inverse $Q$-learning methods.
\end{remark}

\section{DDP-Based Nonlinear Data-Driven IRL}\label{sec:nonlinear}
In this section, we further consider the inverse problem for optimal control problems with a nonlinear cost:
\begin{equation} \label{eq:oc_nonlinear}
 \begin{aligned}
   \mathop {\min }\limits_{u} {\kern 4pt} & \mathbb{E}\Big[\sum_{k=0}^{T-1} l_k(x_k,u_k;\theta) \Big]\\
   s.t.{\kern 5pt}  &x_{k+1}=Ax_k+Bu_k+w_k,
 \end{aligned}
\end{equation}
where $w_k$ is the white noise, and the nonlinear cost function is parameterized by $\theta\in\mathbb{R}^l$. Note that unlike the well-studied LQR case, nonlinear optimal control is usually formulated with a finite horizon for numerical tractability.

To solve nonlinear optimal control problems, differential dynamic programming  \cite{mayne1966second} has been widely applied as a powerful numerical trajectory optimization method by repeatedly optimizing the locally-quadratic models along a nominal trajectory. The recursions are closely related to Pantoja's step-wise Newton's method and exhibit quadratic convergence. To this end, in the setup of this paper, we assume that the expert solves \eqref{eq:oc_nonlinear} with a DDP algorithm. Then the goal of model-free IRL is to learn the parameter $\theta$ from the observed optimal demonstrations $\lbrace x_k^*,u_k^* \rbrace_{k=0}^T$ of a DDP solver, where prior knowledge on the system matrices $(A,B)$ is also unavailable and is approximated by online data collection along \eqref{eq:sys} under PE conditions.

\subsection{Preliminary of DDP}
For stochastic systems with additive process noise, in each iteration of the stochastic DDP, the update rule for the optimal control increment remains the same as that in the deterministic case~\cite{Theodorou2010}.
In each iteration, a backward pass is first performed along the current trajectory to generate a new optimal control sequence, which is then followed by a forward pass to update the nominal trajectory. Before moving on, the basic idea and related notations of DDP are briefly introduced, based on which a data-driven version will be further developed in the next part.

Firstly, regarding the backward pass, let $\lbrace \bar{x}_k, \bar{u}_k \rbrace_{k=0}^T$ denote the current nominal trajectory. The Bellman equation for \eqref{eq:oc_nonlinear} is given by
\begin{equation}\label{eq:Bellman_nonlinear}   
\begin{aligned}
    V_k(x_k)&=\min_{u_k}Q_k(x_k,u_k)\\
    &=\min_{u_k} \lbrace l_k(x_k,u_k,\theta)+V_{k+1}(x_{k+1}) \rbrace.
\end{aligned}
\end{equation}
Considering the 2nd-order Taylor series approximation of \eqref{eq:Bellman_nonlinear} then results in the locally-quadratic optimization model as 
\begin{equation}
    V_k(x_k)=\min_{\delta u_k} \lbrace \delta Q_k(\delta x_k, \delta u_k)+Q_k( \bar{x}_k, \bar{u}_k)) \rbrace,
\end{equation}
where $\delta x$ and $\delta u$ are the variations of the state and control input respectively and 
\begin{equation}\label{eq:delta_Q}
\begin{aligned}
    &\delta Q_k(\delta x_k, \delta u_k)\\
    \approx& \begin{bmatrix}
        Q_{x,k}\\
        Q_{u,k}
    \end{bmatrix}^{\top}\begin{bmatrix}
        \delta x_k\\ \delta u_k
    \end{bmatrix}+ \frac{1}{2}\begin{bmatrix}
        \delta x_k\\ \delta u_k
    \end{bmatrix}^{\top}\begin{bmatrix}
        Q_{xx,k} & Q_{xu,k}\\
        Q_{ux,k} & Q_{uu,k}
    \end{bmatrix}\begin{bmatrix}
        \delta x_k\\ \delta u_k
    \end{bmatrix},
    \end{aligned}
\end{equation}
with
\begin{equation}\label{Q_diff}
    \begin{aligned}
        & Q_{x,k}(\bar{x}_k, \bar{u}_k)=l_{x,k}(\bar{x}_k, \bar{u}_k;\theta)+A^{\top}V_{x,k+1}(\bar{x}_{k+1}),\\
        & Q_{u,k}(\bar{x}_k, \bar{u}_k)=l_{u,k}(\bar{x}_k, \bar{u}_k;\theta)+B^{\top}V_{x,k+1}(\bar{x}_{k+1}),\\
        & Q_{xx,k}(\bar{x}_k, \bar{u}_k)=l_{xx,k}(\bar{x}_k, \bar{u}_k;\theta)+A^{\top}V_{xx,k+1}(\bar{x}_{k+1})A,\\
        & Q_{uu,k}(\bar{x}_k, \bar{u}_k)=l_{uu,k}(\bar{x}_k, \bar{u}_k;\theta)+B^{\top}V_{xx,k+1}(\bar{x}_{k+1})B,\\
        & Q_{ux,k}(\bar{x}_k, \bar{u}_k)=l_{ux,k}(\bar{x}_k, \bar{u}_k;\theta)+B^{\top}V_{xx,k+1}(\bar{x}_{k+1})A,
    \end{aligned}
\end{equation}
and the boundary condition $V_{x,T}=V_{xx,T}=0$. All the partial derivatives in \eqref{eq:delta_Q} are evaluated along the nominal trajectory $\lbrace \bar{x}_k, \bar{u}_k \rbrace_{k=0}^T$ as in \eqref{Q_diff}  and the arguments are omitted in this section for conciseness. In addition, for any function we denote $(\cdot)_{x,k}$ and $(\cdot)_{u,k}$ as its partial derivative at time $k$ with respect to $x$ and $u$, respectively. Then the optimal control increment is derived as 
\begin{equation}\label{eq:delta_u}
\begin{aligned}
    \delta u_k^*&=\mathop{\arg\min}_{\delta u_k} ~\delta Q_k(\delta x_k, \delta u_k)\\
    &=-Q_{uu,k}^{-1}Q_{ux,k}{\delta x_k}-Q_{uu,k}^{-1}Q_{u,k}\\
    &:=K_k{\delta x_k}+k_k
\end{aligned}
\end{equation}

Substituting \eqref{eq:delta_u} back into \eqref{eq:delta_Q} further gives the updating rule of the value function at $\bar{x}_k, \bar{u}_k$ as:
\begin{equation}\label{eq:V_diff}
    \begin{aligned}
        &V_{x,k}=Q_{x,k}-Q_{ux,k}^{\top}Q_{uu,k}^{-1}Q_{u,k},\\
        &V_{xx,k}=Q_{xx,k}-Q_{ux,k}^{\top}Q_{uu,k}^{-1}Q_{ux,k}.
    \end{aligned}
\end{equation}

Therefore, starting at $t=N$ and solving the locally-quadratic approximation of the Bellman equation backwards in time via (\ref{Q_diff})-(\ref{eq:V_diff}), the optimal control sequence $\lbrace \delta u_N^*, \delta u_{N-1}^*, \cdots, \delta u_0^* \rbrace$ is obtained at the current nominal trajectory $\lbrace \bar{x}_k, \bar{u}_k \rbrace_{k=0}^T$. Then the nominal trajectory can be updated via a forward pass as
\begin{equation}
    \delta x_{k+1} = A \delta x_{k} + B\delta u_{k},~~\delta x_{0}=0,
\end{equation}
and $x_{k}=\bar{x}_{k}+\delta{x}_{k}$ for $k-0,1,\cdots,T$.

By repeatedly applying the backward and forward passes, the nominal trajectory will finally converge to the optimal state trajectory to \eqref{eq:oc_nonlinear}.

\subsection{DDP-Based Model-Free Nonlinear IRL}
In each iteration of DDP, a quadratic optimization problem is solved such that the optimal control increment is obtained in a manner similar to the LQR problem. However, developing data-driven characterization of DDP iterations differs from the model-free formulation of the ARE in Section \ref{sec:LQR} in two aspects: (i) the quadratic cost in \eqref{eq:delta_Q} also involves a cross-penalty of $(\delta x_k,\delta u_k)$ and an affine term; (ii) a time-varying optimal feedback law is involved as in \eqref{eq:delta_u}.

To begin with, a data-driven reformulation of the backward pass is studied. To obtain the partial derivatives $V_{(\cdot)}$ iteratively along the current nominal trajectory, we denote 
\begin{subequations}\label{eq:Hk_ck}
     \begin{align}
        H_k&:=\begin{bmatrix}
            l_{xx,k}+A^{\top}V_{xx,k+1}A & A^{\top}V_{xx,k+1}B+l_{xu,k}\\
            B^{\top}V_{xx,k+1}A+l_{ux,k}   & l_{uu,k}+B^{\top}V_{xx,k+1}B
        \end{bmatrix}\notag\\
        &:=\begin{bmatrix}
        H_{1,k} & H_{2,k}^{\top} \\
        H_{2,k} & H_{3,k}
        \end{bmatrix}, \label{eq:Hk}\\
        c_k&:=\begin{bmatrix}
            l_{x,k}+A^{\top}V_{x,k+1}\\
            l_{u,k}+B^{\top}V_{x,k+1}
        \end{bmatrix}. \label{eq:ck}
     \end{align}
\end{subequations}
     
\begin{proposition}\label{prop:Hk_recur}
    If $\lbrace V_{x,k}, V_{xx,k} \rbrace_{k=0}^T$ is obtained following the recursions in \eqref{Q_diff} and \eqref{eq:V_diff}, then $\lbrace H_{k} \rbrace_{k=0}^{T-1}$ and $\lbrace c_{k} \rbrace_{k=0}^{T-1}$  satisfy the backwards recursions:
\begin{subequations}\label{eq:H_c_iter}
     \begin{align}
        &H_k=\begin{bmatrix}
            l_{xx,k} & l_{xu,k}\\
            l_{ux,k} & l_{uu,k}
        \end{bmatrix}+\Xi_{k+1}^{\top}H_{k+1}\Xi_{k+1},\label{eq:Hk_recursion}\\
        &c_k=\begin{bmatrix} l_{x,k}\\l_{u,k} \end{bmatrix}
        +\begin{bmatrix}
            A & B\\
            K_{k+1}A & K_{k+1}B
        \end{bmatrix}^{\top}c_{k+1}, \label{eq:ck_recursion}
     \end{align}
\end{subequations}
    where $K_{k+1}=-H_{3,k+1}^{-1}H_{2,k+1}$ and
    \begin{equation*}
        \Xi_{k+1}:=\Xi_{k+1}(H_{k+1})=\begin{bmatrix}
            A & B\\
            K_{k+1}A & K_{k+1}B
        \end{bmatrix}.
    \end{equation*}
\end{proposition}
\proof
    Since the backward recursions remain the same as that in deterministic systems, for brevity of notations we let $w_k=0,~\forall k\geq 0$ in the proof. Let $\lbrace u_k {\rbrace}_{k=0}^{T-1}$ be an arbitrary input sequence that is persistently exciting and $\lbrace x_k {\rbrace}_{k=0}^{T-1}$ be the corresponding state trajectory. Denote $z_k=[x_k^{\top}~u_k^{\top}]^{\top}$. Then straightforward computations on \eqref{eq:Hk} yield
    \begin{equation}\label{eq:zHz}
    \begin{aligned}
        &z_i^{\top}H_kz_j\\
        =&z_i^{\top}\begin{bmatrix}
        l_{xx,k} & l_{xu,k}\\
        l_{ux,k} & l_{uu,k}
    \end{bmatrix}z_j+(Ax_{i}+Bu_{i})^{\top}V_{xx,k+1}(Ax_{j}+Bu_{j}).        
    \end{aligned}
    \end{equation}
    Based on \eqref{Q_diff} and \eqref{eq:V_diff}, one can show that
    \begin{equation}\label{eq:DRE_nonlin}
    \begin{aligned}
        V_{xx,k}
        =&l_{xx,k}+(A+BK_k)^{\top}V_{xx,k+1}(A+BK_k)\\
        &+K_k^{\top}l_{uu,k}K_k+K_k^{\top}l_{ux,k}+l_{ux,k}^{\top}K_k,
    \end{aligned}
    \end{equation}
    Similar to \eqref{eq:zHz},  we evaluate \eqref{eq:Hk} at one time instant forward, which together with \eqref{eq:DRE_nonlin} further yields
    \begin{equation*}
        \begin{aligned}
            &\begin{bmatrix}
                x_{i+1} \\ K_{k+1}x_{i+1}
            \end{bmatrix}^{\top}H_{k+1}\begin{bmatrix}
                x_{j+1} \\ K_{k+1}x_{j+1}
            \end{bmatrix}\\
            =&\begin{bmatrix}
                x_{i+1} \\ K_{k+1}x_{i+1}
            \end{bmatrix}^{\top}\begin{bmatrix}
        l_{xx,k+1} & l_{xu,k+1}\\
        l_{ux,k+1} & l_{uu,k+1}
    \end{bmatrix}\begin{bmatrix}
                x_{j+1} \\ K_{k+1}x_{j+1}
            \end{bmatrix}\\
            &+x_{i+1}^{\top}(A+BK_{k+1})^{\top}V_{xx,k+2}(A+BK_{k+1})x_{j+1}\\
    =&x_{i+1}^{\top}V_{xx,k+1}x_{j+1}.
        \end{aligned}
    \end{equation*}
    Therefore, \eqref{eq:zHz} can be rewritten as 
    \begin{equation*}\label{eq:Hk_recur_data}
    \begin{aligned}
        &z_i^{\top}H_kz_j\\
        =&z_i^{\top}\begin{bmatrix}
        l_{xx,k} & l_{xu,k}\\
        l_{ux,k} & l_{uu,k}
        \end{bmatrix}z_j+\begin{bmatrix}
                x_{i+1} \\ K_{k+1}x_{i+1}
            \end{bmatrix}^{\top}H_{k+1}\begin{bmatrix}
                x_{j+1} \\ K_{k+1}x_{j+1}
            \end{bmatrix}.
    \end{aligned}
    \end{equation*}
    Let $Z=[z_0~z_1~\cdots~ z_{n+m-1}]$. Since $x_{l+1}=Ax_l+Bu_l$ with $l=i,j$, the above equation further implies
    \begin{equation}
    Z^{\top}(H_k-\begin{bmatrix}
        l_{xx,k} & l_{xu,k}\\
        l_{ux,k} & l_{uu,k}
        \end{bmatrix}-\Xi_{k+1}^{\top} H_{k+1}\Xi_{k+1})Z=0.
    \end{equation}
    Under the assumption of PE inputs, we know that $Z$ has full row rank. Hence, \eqref{eq:Hk_recursion} follows. In addition, the recursions of $c_k$ in \eqref{eq:ck_recursion} can be obtained directly by substituting \eqref{Q_diff} and \eqref{eq:V_diff}.
    Therefore, the proof is completed.
\endproof

We next show that \eqref{eq:H_c_iter} indeed provides an equivalent formula for the backwards pass in \eqref{Q_diff} and \eqref{eq:V_diff}.
\begin{theorem} \label{thm:Hk_DRE}
    Given the boundary condition
\begin{equation}\label{eq:Hc_bound}
    H_{T-1}=\begin{bmatrix}
        l_{xx,T-1} & l_{xu,T-1}\\
        l_{ux,T-1} & l_{uu,T-1}
    \end{bmatrix},~c_{T-1}=\begin{bmatrix} l_{x,T-1}\\l_{u,T-1} \end{bmatrix},
\end{equation}
the sequence $\lbrace H_k, c_k \rbrace_{k=0}^{T-1}$ satisfies the backwards recursion in \eqref{eq:H_c_iter} if and only if \eqref{eq:Hk_ck} holds, where $V_{x,k+1}$ and $V_{xx,k+1}$ are defined by \eqref{Q_diff} and \eqref{eq:V_diff}.
\end{theorem}
\proof
    The sufficiency is shown in Proposition \ref{prop:Hk_recur}. Next, we show the necessity by induction. Suppose $\lbrace H_k, c_k \rbrace_{k=0}^{T-1}$ is solved from the backwards recursion in \eqref{eq:H_c_iter}. It is clear that \eqref{eq:Hk_ck} holds at $k=T-1$. Assume that $H_{k+1}$ satisfies \eqref{eq:Hk}, namely,
    \begin{equation*}\label{eq:Hk1_Q}
        H_{k+1}=\begin{bmatrix}
        Q_{xx,k+1} & Q_{xu,k+1}\\
        Q_{ux,k+1} & Q_{uu,k+1}
    \end{bmatrix}.
    \end{equation*}
    Substituting $K_{k+1}=-H_{3,k+1}^{-1}H_{2,k+1}$, then straightforward computations yield 
    \begin{equation}
        \begin{aligned}
            &\begin{bmatrix}
                A & B\\
                K_{k+1}A & K_{k+1}B
            \end{bmatrix}^{\top}H_{k+1}\begin{bmatrix}
                A & B\\
                K_{k+1}A & K_{k+1}B
            \end{bmatrix}\\
            =&\begin{bmatrix}
                A^{\top}\\ B^{\top}
            \end{bmatrix}(H_{1,k+1}-H_{2,k+1}^{\top}H_{3,k+1}^{-1}H_{2,k+1})\begin{bmatrix}
                A & B
            \end{bmatrix}\\
            =&\begin{bmatrix}
                A^{\top}\\ B^{\top}
            \end{bmatrix}V_{xx,k+1}\begin{bmatrix}
                A & B
            \end{bmatrix}.
        \end{aligned}
    \end{equation}
    Hence, the right-hand side of \eqref{eq:Hk_recursion} is equal to that of \eqref{eq:Hk}, which means that $H_{k}$ also satisfies \eqref{eq:Hk}. Furthermore, substituting \eqref{eq:Hk} into \eqref{eq:ck_recursion} directly results in its equivalence to \eqref{eq:ck}.
    Therefore, the necessity is proved.
\endproof

Similar to the data-driven  ILQR in Section \ref{sec:LQR}, the backward update rule of $H_k$ can also be implemented based on system data $(X_{0,t},U_{0,t},X_{1,t})$ without using system matrices $(A,B)$ as:
\begin{equation}\label{eq:Hk_sol}
    H_k=\begin{bmatrix}
        l_{xx,k} & l_{xu,k}\\
        l_{ux,k} & l_{uu,k}
        \end{bmatrix}+\hat{\Xi}_{k+1,t}^{\top}H_{k+1}\hat{\Xi}_{k+1,t}.
\end{equation}
where $\hat{\Xi}_{k+1,t}=\begin{bmatrix}
            I_n& {K}_{k+1}^{\top}
        \end{bmatrix}^{\top}\bar{X}_{1,t}\Lambda_t^{-1}$.

In addition, with given $c_{k+1}$ and $K_{k+1}$, $c_k$ can be obtained by \eqref{eq:ck_recursion} in a model-free manner as 
\begin{equation}\label{eq:c_data}
c_k=\begin{bmatrix} l_{x,k}\\l_{u,k} \end{bmatrix}
        +\hat{\Xi}_{k+1,t}^{\top}c_{k+1}.
\end{equation}
\begin{remark}
    In real implementations of the DDP solver, line-search methods or regularization terms can be involved to guarantee the positive-definiteness of $Q_{uu,k}$. Such techniques have been well-studied in the literature, whose details are omitted in this paper. Here we only assume that when the DDP terminates, the matrices $\lbrace Q_{uu,k} \rbrace_{k=0}^{T-1}$ are all positive definite along the optimal solution (i.e., the observed demonstrations). Together with some mild condition on the SNR, such assumption ensures the well-posedness of the data-driven parameterization in \eqref{eq:H_c_iter}.
\end{remark}

Now we are ready to formulate a DDP-based model-free IRL problem for \eqref{eq:oc_nonlinear}. Let $\theta^*$ denote the true parameter of the cost that generates the observed demonstrations $\lbrace x_k^*,u_k^* \rbrace_{k=0}^T$. In the sequel, we define
$\bar{l}_{k}(\theta):=l_k(x_k^*,u_k^*;\theta)$ as the evaluation of the cost function at some candidate parameter $\theta$ along the demonstrations $\lbrace x_k^*,u_k^* \rbrace_{k=0}^T$, and the same notations also apply to the partial derivatives of $l_k(\cdot)$. Recall that in the DDP solver we have $\bar{Q}_{u,k}(\theta^*):=Q_{u,k}({x}_k^*,{u}_k^*;\theta^*)=0$, which implies $\norm{\bar{Q}_{u,k}(\theta)}$ can be used as the residual error for the inverse problem. Therefore, based on the data parameterization of the DDP recursions \eqref{eq:Hk_sol} and \eqref{eq:c_data}, a data-driven formulation of the inverse optimal control problem to \eqref{eq:oc_nonlinear} is given by
\begin{equation} \label{eq:nonlinear_IOC}
 \begin{aligned}
   \mathop {\min }\limits_{\theta\in \mathbb{R}^l} {\kern 4pt} & L(\theta)= \sum_{k=0}^{T-1} \norm{\begin{bmatrix} 0& I_m\end{bmatrix} c_k}^2 \\
   s.t.{\kern 5pt}  & H_k=\Gamma_k(\theta)+\hat{\Xi}_{k+1,t}^{\top}H_{k+1}\hat{\Xi}_{k+1,t},\\
        &(k=0,1,\cdots,T-2)\\
   &c_k=\eta_k(\theta)
        +\hat{\Xi}_{k+1,t}^{\top}c_{k+1},~(k=0,1,\cdots,T-2)\\
   &H_{T-1}=\Gamma_{T-1}(\theta),~c_{T-1}=\eta_{T-1}(\theta),
 \end{aligned}
\end{equation}
where $\hat{\Xi}_{k+1,t}$ is a function of $H_{k+1}$ as defined in \eqref{eq:Hk_sol} and 
\begin{equation*}
    \Gamma_k(\theta)=\begin{bmatrix}
        \bar{l}_{xx,k}(\theta) & \bar{l}_{xu,k}(\theta)\\
        \bar{l}_{ux,k}(\theta) & \bar{l}_{uu,k}(\theta)
        \end{bmatrix},~\eta_k(\theta)=\begin{bmatrix} \bar{l}_{x,k}(\theta)\\\bar{l}_{u,k}(\theta) \end{bmatrix}.
\end{equation*}

For any $\theta\in\mathbb{R}^l$, $\lbrace H_k,c_k \rbrace_{k=0}^{T-1}$ are uniquely determined as a function of $\theta$ using the backward recursions in the constraints of \eqref{eq:nonlinear_IOC}. Therefore, we only consider the residual cost $L(\theta)$ in the above formulation as an implicit function of $\theta$. With its gradient $\nabla L(\theta)$ given as follows in Proposition \ref{prop:grad_L}, one can easily develop convergent algorithms to recover a best-fit estimate of $\theta$ with the offline data $(X_{0,t},U_{0,t},X_{1,t})$ in a model-free manner. 
\begin{proposition}\label{prop:grad_L}
    For each candidate $\theta=[\theta_1,\theta_2,\cdots,\theta_l]^{\top}$, the gradient of $L(\theta)$ in \eqref{eq:nonlinear_IOC}, namely,
    $$ \nabla L(\theta)= 2\sum_{k=0}^{T-1}c_k^{\top}\begin{bmatrix} 0\\ I_m\end{bmatrix}\begin{bmatrix} 0& I_m\end{bmatrix} \frac{\partial c_k}{\partial \theta}, $$
    can be obtained with backward recursions in Algorithm \ref{alg:grad_L}.
\end{proposition}
\proof
    Differentiating the first two constraints in \eqref{eq:nonlinear_IOC} gives the backward recursions for $\frac{\partial H_k}{\partial \theta_i}$ and $\frac{\partial c_k}{\partial \theta_i}$ as
    \begin{equation}\label{eq:diff_H_c}
    \begin{aligned}
        & \frac{\partial H_k}{\partial \theta_i}=\frac{\partial }{\partial \theta_i}\Gamma_k(\theta)+\frac{\partial }{\partial \theta_i}(\hat{\Xi}_{k+1,t}^{\top}H_{k+1}\hat{\Xi}_{k+1,t}),\\
   &\frac{\partial c_k}{\partial \theta_i}=\frac{\partial }{\partial \theta_i}\eta_k(\theta)
        +\frac{\partial }{\partial \theta_i}(\hat{\Xi}_{k+1,t}^{\top}c_{k+1}),
    \end{aligned}        
    \end{equation}
    where $i=1,2,\cdots,l$ and $k=0,1,\cdots,T-2$, with the terminal boundary given by 
    $$\frac{\partial H_{T-1}}{\partial \theta_i}=\frac{\partial }{\partial \theta_i}\Gamma_{T-1}(\theta),~\frac{\partial c_{T-1}}{\partial \theta_i}=\frac{\partial }{\partial \theta_i}\eta_{T-1}(\theta).$$
    Using the partitions of $H_{k+1}$, we have
    \begin{align*}
        &\hat{\Xi}_{k+1,t}^{\top}H_{k+1}\hat{\Xi}_{k+1,t}\\
        =& \Lambda_t^{-1}\bar{X}_{1,t}^{\top}(H_{1,k+1}-H_{2,k+1}^{\top}H_{3,k+1}^{-1}H_{2,k+1})\bar{X}_{1,t}\Lambda_t^{-1},
    \end{align*}
    which implies 
    \begin{equation}\label{eq:diff_Hk}
    \begin{aligned}
        & \frac{\partial H_k}{\partial \theta_i}=\frac{\partial }{\partial \theta_i}\Gamma_k(\theta)+\Lambda_t^{-1}\bar{X}_{1,t}^{\top}N_{i,k+1}\bar{X}_{1,t}\Lambda_t^{-1},
    \end{aligned}
    \end{equation}
    with
    \begin{equation}  \label{eq:Nk}
    \begin{aligned}
     N_{i,k+1}:=& \frac{\partial H_{1,k+1}}{\partial \theta_i}-2\overline{H_{2,k+1}^{\top}H_{3,k+1}^{-1}\frac{\partial H_{2,k+1}}{\partial \theta_i}}\\
     &+H_{2,k+1}^{\top}H_{3,k+1}^{-1}\frac{\partial H_{3,k+1}}{\partial \theta_i} H_{3,k+1}^{-1}H_{2,k+1}.
    \end{aligned}
    \end{equation}
    Similarly, substituting the expression of $\hat{\Xi}_{k+1,t}$ in \eqref{eq:diff_H_c} further simplifies $\frac{\partial c_k}{\partial \theta_i}$ by
    \begin{equation}\label{eq:diff_ck}
    \begin{aligned}
        \frac{\partial c_k}{\partial \theta_i}=&\frac{\partial \eta_k(\theta)}{\partial \theta_i}
        +\Lambda_t^{-1}\bar{X}_{1,t}^{\top}\begin{bmatrix}
            0_{n\times n} & M_{i,k+1}
        \end{bmatrix}c_{k+1}\\        
        &+\Lambda_t^{-1}\bar{X}_{1,t}^{\top}\begin{bmatrix}
        I & -H_{2,k+1}^{\top}H_{3,k+1}^{-1}
        \end{bmatrix}\frac{\partial c_{k+1}}{\partial \theta_i},
    \end{aligned}
    \end{equation}
    with
    \begin{equation}\label{eq:Mk}
        M_{i,k+1}:=H_{2,k+1}^{\top}H_{3,k+1}^{-1}\frac{\partial H_{3,k+1}}{\partial \theta_i} H_{3,k+1}^{-1}-\frac{\partial H_{2,k+1}^{\top}}{\partial \theta_i} H_{3,k+1}^{-1}.
    \end{equation}
    Summarizing the above results, one can derive Algorithm \ref{alg:grad_L}.
\endproof
\begin{algorithm}[t]
\caption{Gradient solver for $L(\theta)$.}\label{alg:grad_L}
\begin{algorithmic}
    \REQUIRE {Candidate $\theta$, system data $\Lambda_t$ and $\bar{X}_{1,t}$.}
    \ENSURE {$\nabla L(\theta)$ }
    \STATE Initialize: $H_{T-1}=\Gamma_{T-1}(\theta),~c_{T-1}=\eta_{T-1}(\theta)$.
    \FOR{$k=T-2,\cdots,0$} 
        \FOR{$i=1,\cdots,l$}
        \STATE Obtain $N_{i,k+1}$ and $M_{i,k+1}$ using \eqref{eq:Nk} and \eqref{eq:Mk}.
        \STATE Compute $\frac{\partial H_k}{\partial \theta_i}$ and $\frac{\partial c_k}{\partial \theta_i}$ via \eqref{eq:diff_Hk} and \eqref{eq:diff_ck}.      
        \ENDFOR
        \STATE Update $H_k$ and $c_k$ by
        \begin{align*}
            & H_k=\Gamma_k(\theta)+\hat{\Xi}_{k+1,t}^{\top}H_{k+1}\hat{\Xi}_{k+1,t},\\
            &c_k=\eta_k(\theta)
            +\hat{\Xi}_{k+1,t}^{\top}c_{k+1}.
        \end{align*}
    \ENDFOR
    \RETURN 
    $ \nabla L(\theta)= 2\sum_{k=0}^{T-1}c_k^{\top}\diag(0_{n\times n}, I_m)\Big[\frac{\partial c_k}{\partial \theta_1}~\cdots~\frac{\partial c_k}{\partial \theta_l} \Big].$ 
\end{algorithmic}
\end{algorithm}    

It is noteworthy that when system data is collected in an online manner, the dimensions of variables in \eqref{eq:nonlinear_IOC} are invariant as $t$ increases. Then similar to the LQR case in Section \ref{sec:LQR}, one may also develop online iterations of $\theta$ based on Proposition \ref{prop:grad_L}, where various techniques in online optimization can be borrowed. However, it is impossible to derive the convergence properties of the optimization algorithm for general nonlinear problems, which indeed depend on the specific structure and properties of $l(x_k,u_k;\theta)$. 
Since the main focus of this section is to provide a general and scalable data-driven framework for nonlinear IRL that enables online data collection, detailed convergence analysis for specific cost structures is left for our future work due to the page limitations. Performances of the nonlinear IRL are demonstrated via numerical simulations.

\section{Numerical Simulation}\label{sec:simulation}
\subsection{Convergence of Adaptive Learning in ILQR}
In this section, we apply Algorithm \ref{alg:Newton_data} for adaptive IRL in the LQR with online PE data. The following system parameters are randomly generated to obtain the expert demonstration $\bar K$
\begin{equation}\label{equ:simulation_matrix}
\begin{aligned}
&A= \begin{bmatrix}1 & 0.1 \\ -0.2946 & 1.1838   \end{bmatrix}, \, B=\begin{bmatrix} -0.1098 \\ -2.2375\end{bmatrix},    \\
&Q = \begin{bmatrix} 0.5029 & -0.2047 \\ -0.2047 & 0.3384\end{bmatrix}, \, x_{0} = \begin{bmatrix}2 \\ 1    \end{bmatrix}.
\end{aligned}
\end{equation}

\begin{figure}[!t]
    \centering
    \subfloat[Adaptive ILQR]{%
        \includegraphics[height=4.2cm]{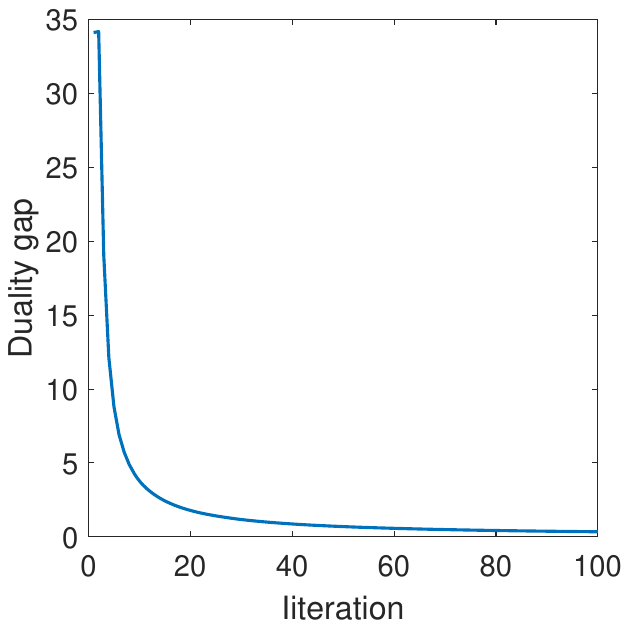}}
    \subfloat[DDP-based IRL]{%
        \includegraphics[height=4.2cm]{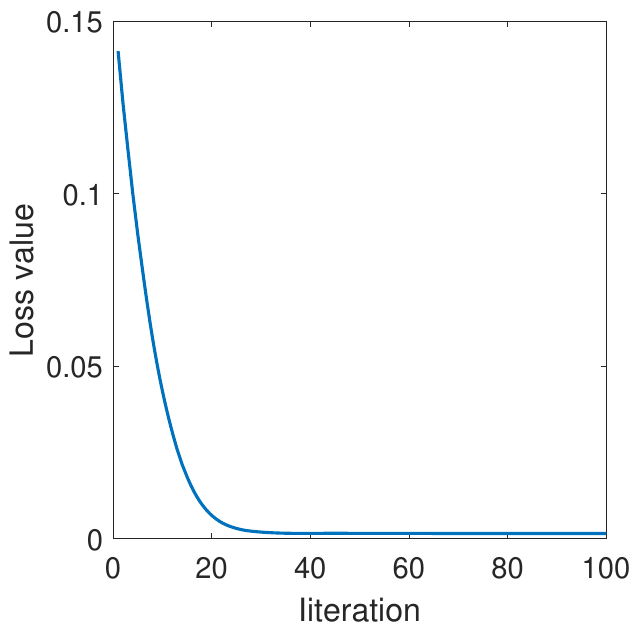}}       
    \caption{Convergence of online IRL.}
    \label{fig:convergence} 
\end{figure}

Without prior knowledge of $(A,B)$, online system data is generated under the PE condition to adaptively update the estimate of $Q$ via the data-driven IRL algorithm.
The initial of $Q$ is chosen such that $\delta_{0} < \frac{1}{2}-h_7$. We choose a large initial duality gap to illustrate that our online algorithm possesses a quite large region of convergence. With $t_0=20$ and applying the online algorithm, convergence of the duality gap is shown in Fig. \ref{fig:convergence}(a), thereby demonstrating the convergence of Algorithm \ref{alg:Newton_data}. Moreover, although the data matrix $\hat{\Phi}_{t}$ is updated at each step as new data is incorporated, the convergence of Algorithm \ref{alg:Newton_data} remains smooth throughout the process.

\begin{figure}[t]
    \centering
    \includegraphics[width=0.75\linewidth]{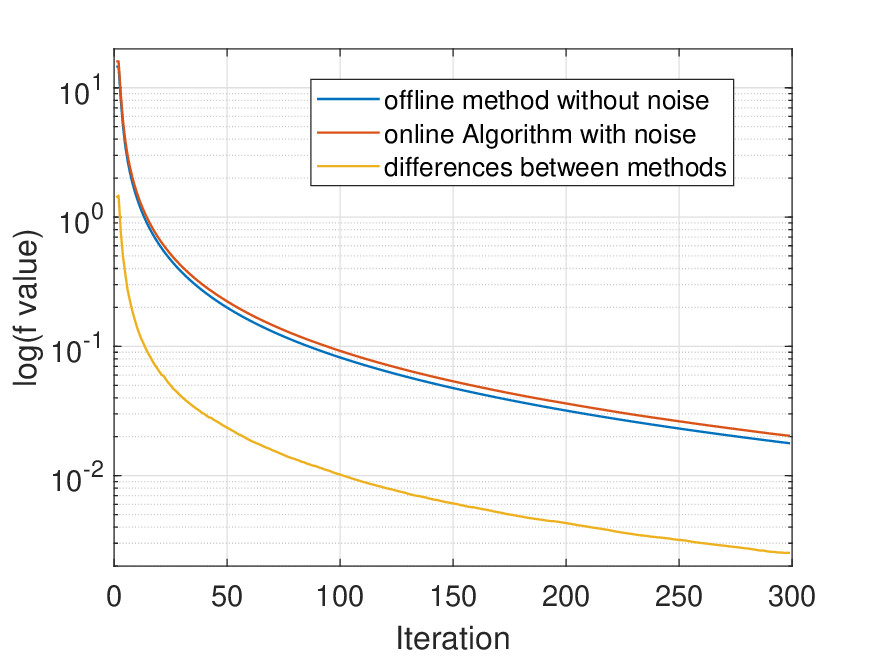}
    \caption{Comparison of adaptive IRL with model-based IRL.}
    \label{fig:compare_with_offline}
\end{figure}

Furthermore, applying the model-based IRL with the true system matrices $(A,B)$ as the ground truth, the accuracy of the online cost estimate is also tested by evaluating the value of $f(Q)$ in \eqref{eq:LS_Q_K_given}, which in essence characterizes the residual error between $\bar{K}$ and the reconstructed optimal policy based on the true system model.
We conduct 300 iterations for both the online data-driven algorithm and the model-based method. 
For clarity of presentation, $f(Q)$ is plotted on a logarithmic scale. As shown in Figure \ref{fig:compare_with_offline}, the discrepancy between Algorithm \ref{alg:Newton_data} and the model-based method diminishes with increasing iterations, which indicates that our proposed algorithm achieves a comparable level of performance despite the presence of noise.

\subsection{DDP-Based IRL in Human-Aware Robot Navigation}
\begin{figure}[t]
    \centering
    \includegraphics[width=0.66\linewidth]{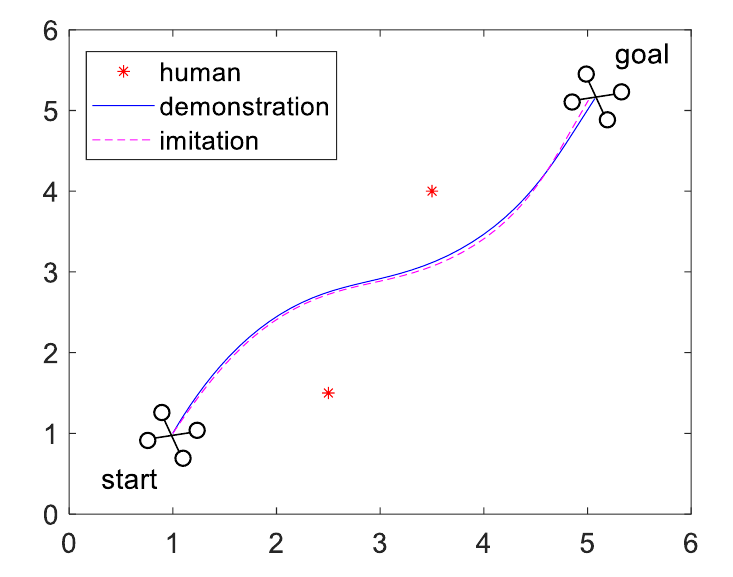}
    \caption{Comparison of optimal waypoints.}
    \label{fig:pathplan}
\end{figure}

We consider the human-aware autonomous navigation task of a drone, whose goal is to navigate in an unknown environment by following a nominal trajectory while avoiding some humans. The navigation tasks of UAVs are typically achieved in a cascaded framework, comprising a high-level path planner that optimizes waypoints and a low-level actuation module (such as a PID controller) to track the commands. Let $p_k, v_k, u_k\in\mathbb{R}^2$ denote the position, velocity and acceleration of the robot, respectively. The path planner generates the waypoints by minimizing the instantaneous cost:
$$l_k=\theta_u\norm{u_k}^2+\theta_w\norm{p_k-p_{w,k}}^2+\sum_{i=1}^M\theta_1e^{(-\theta_2\norm{p_k-p_{h,i}}^2)},$$
where $\{p_{h,i}\}_{i=1}^M$ denote the positions of $M$ humans, $\{p_{w,k}\}_{k=0}^T$ denote the nominal waypoints to follow in the absence of humans, and $\theta_{(\cdot)}$ are positive weights.
The above squared-exponential distance function has been widely used to describe the personal space in social navigation \cite{KollmitzIROS}. Let $\theta_u=0.01$ be fixed. For normalization, the state penalties satisfy $\theta_w+M\theta_1=1$. Through a DDP solver, an example of the optimal path is shown by the blue line in Fig. \ref{fig:pathplan}. 

Then the online IRL algorithm is applied to reconstruct $\theta=[\theta_1,\theta_2,\theta_w]$ from demonstrations. We also choose $t_0=20$ and a constant step-size is applied. The convergence of the loss value is shown in Fig. \ref{fig:convergence}(b).
To illustrate the learning accuracy, the forward DDP is solved with the recovered cost parameters and the regenerated optimal path is illustrated by the dotted line in Fig. \ref{fig:pathplan}, which imitates the demonstrations quite well. 

\section{Conclusion}\label{sec:conclusion}
In this paper, we addressed the online IRL problems for reconstructing unknown cost functions from observed optimal behavior in a model-free manner. A novel adaptive algorithm was developed using completely off-policy system data, requiring only a mild PE condition. By employing full NT-step primal-dual interior-point iterations, the algorithm directly updates the cost parameter without a complete RL procedure as in conventional bi-level methods, and its non-asymptotic convergence with a sublinear rate was rigorously established. Numerical simulations show that the proposed algorithm can achieve performance comparable to offline noiseless benchmarks. Furthermore, the method was generalized to nonlinear IRL based on a DDP framework, which enables direct computation of the loss gradient in a model-free manner.

\appendix
\section{Appendix}   
\subsection{Proof of Lemma \ref{lem:Phi_t_stable}}\label{append:Phi_t_stable}
By Lemma 2 in \cite{lopez2023efficient}, it holds that $\hat\Phi_t$ is stable if and only if $\hat{A}_t-\hat{B}_t\bar{K}$ is stable. Since the system model can be rewritten as $X_{1,t}=AX_{0,t}+BU_{0,t}+W_{0,t}$, right-multiplying both sides by $\frac{1}{t}Z_{0,t}^{\top}$ yields
    \begin{equation*}
        \bar{X}_{1,t}=\begin{bmatrix}
            A & B
        \end{bmatrix}\Lambda_t+\bar{W}_{0,t},
    \end{equation*}
    which under the PE condition further implies
    \begin{equation}\label{eq:AB}
        \begin{bmatrix}
            A & B
        \end{bmatrix}=(\bar{X}_{1,t}-\bar{W}_{0,t})\Lambda_t^{-1}.
    \end{equation}
    Since $[\hat{A}_t~\hat{B}_t]=X_{1,t}Z_{0,t}^{\dagger}=\bar{X}_{1,t}\Lambda_t^{-1}$, we have
    \begin{equation}\label{eq:Acl_err}
        \begin{aligned}
            \hat{A}_t-\hat{B}_t\bar{K}=(A-B\bar{K})+\bar{W}_{0,t}\Lambda_t^{-1}\begin{bmatrix}
                I_n & -\bar{K}^{\top}
            \end{bmatrix}^{\top}.
        \end{aligned}
    \end{equation}
By the robust fundamental lemma \cite{ZhaoData2025},  Assumption \ref{assump:PE} ensures a lower bound for the singular values of $Z_{0,t}$ by $\underline{\sigma}(Z_{0,t})\geq\sqrt{t}\gamma\zeta/2$, which implies
    \begin{equation}\label{eq:WZ_bnd}
        \norm{\bar{W}_{0,t}\Lambda_t^{-1}}_2=\norm{{W}_{0,t}{Z}_{0,t}^{\dagger}}_2\leq\frac{\norm{{W}_{0,t}}_2}{\underline{\sigma}(Z_{0,t})}\leq \frac{2\kappa}{\gamma\zeta}.
    \end{equation}     
    Under the assumption of \eqref{eq:SNR_bnd}, \eqref{eq:Acl_err} further yields
    \begin{align*}
        &\norm{(\hat{A}_t-\hat{B}_t\bar{K})-(A-B\bar{K})}_2\\
        =&\norm{\bar{W}_{0,t}\Lambda_t^{-1}\begin{bmatrix}
                I_n & -\bar{K}^{\top}
            \end{bmatrix}^{\top}}_2\\
            \leq& \norm{\bar{W}_{0,t}\Lambda_t^{-1}}_2\norm{\begin{bmatrix}
                I_n & -\bar{K}^{\top}
            \end{bmatrix}}_2\\
            \leq& \frac{2\kappa}{\gamma \zeta}(\norm{\bar{K}}_2+1) \leq (\sqrt{\bar{\sigma}(\Sigma_K-I_n)\bar{\sigma}(\Sigma_K)}+\bar{\sigma}(\Sigma_K))^{-1}.
    \end{align*}
    Following the analysis of robust stability in \cite{sezer1988robust}, we can guarantee that $\hat{A}_t-\hat{B}_t\bar{K}$ is also stable under \eqref{eq:SNR_bnd}.

\subsection{Proof of Proposition \ref{prop:convF_grad}}\label{append:convF_grad}
    In order to show that $f_t(Q)$ is convex, we consider an arbitrary direction $\tilde{Q} \in \mathbb{S}^n$ and define $\tilde H$ as the unique solution that solves
    $ \tilde{H} = \diag (\tilde{Q},0_{m\times m}) +\hat{\Phi}_t^{\top}\tilde{H}\hat{\Phi}_t. $
    Then straightforward computations yield that 
    $$f_t(Q+\alpha\tilde{Q})=\frac{1}{2}	\langle H+\alpha\tilde{H},\Psi(H+\alpha\tilde{H}) \rangle,$$
    and 
    $$\frac{d}{d\alpha}f_t(Q+\alpha\tilde{Q})=\langle \tilde{H},\Psi(H+\alpha\tilde{H}) \rangle.$$
    Hence, the Hessian of $f_t(Q+\alpha\tilde{Q})$ at $Q$ along the direction $\tilde{Q}$ is given by
    \begin{align*}
        \frac{d^2}{d\alpha^2}f_t(Q+\alpha\tilde{Q})\vert_{\alpha=0}&=\langle \tilde{H},\Psi(\tilde{H}) \rangle\\
        &=\norm{\begin{bmatrix}
        0_{m\times n} & I_m
    \end{bmatrix} \tilde{H}\begin{bmatrix}
        I_n \\ -\bar{K}
    \end{bmatrix}}^2\geq 0,
    \end{align*}
    which holds for any $\tilde{Q} \in \mathbb{S}^n$. Therefore, $f_t(Q)$ is convex.

    Next, the gradient of $f(Q)$ is computed by analyzing its differential. Computing the differential of \eqref{eq:H_lyap_t} yields
    \begin{equation*}
    dH=\diag(dQ,0_{m\times m})+\hat{\Phi}_t^{\top}dH\hat{\Phi}_t. 
    \end{equation*}
    Since $\hat{\Phi}_t$ is stable, the unique solution of $dH$ to the above Lyapunov equation always exists, which is given by 
    $$ dH=\sum_{i=0}^{\infty} (\hat{\Phi}_t^{\top})^i\diag(dQ,0_{m\times m}) \hat{\Phi}_t^i. $$
    Similarly, for any $H(Q)$, the unique solution of $Z$ that solves \eqref{eq:Z_syl} is given by $Z=\sum_{i=0}^{\infty}\hat{\Phi}_t^{i}\overline{\Psi(H(Q))}(\hat{\Phi}_{t}^{\top})^i$.   
    Then the differential $df_t$ can be computed as follows:
    \begin{equation*}
        \begin{aligned}
            df_t(Q)&=\frac{1}{2}\langle dH(Q),\Psi(H(Q)) \rangle+\frac{1}{2}\langle H(Q),d\Psi(H(Q)) \rangle\\
            &=\tr (\overline{\Psi(H(Q))}dH)\\
            &=\tr (\sum_{i=0}^{\infty}\hat{\Phi}_t^i\overline{\Psi(H(Q))}(\hat{\Phi}_t^{\top})^i\diag(dQ,0_{m\times m}))\\
        &=\tr (Z\diag(dQ,0_{m\times m}))=\tr (\begin{bmatrix}
        I_n & 0
        \end{bmatrix}Z\begin{bmatrix}
            I_n \\ 0
        \end{bmatrix}dQ).
        \end{aligned}
    \end{equation*}
    Since $df_t(Q)=\langle \nabla f_t(Q), dQ \rangle$, the expression of $\nabla f_t(Q)$ follows directly as in \eqref{eq:Delta_f} and \eqref{eq:Z_syl}.
    
\subsection{Proof of Proposition \ref{prop:dual}}\label{sec:append_prop_dual}
As for the primal problem in \eqref{eq:Jt}, its Lagrange dual function defined on $S\in\mathbb{S}^n_+$ is given by
\begin{align*}
    \varphi_t(S)&:=\inf_{Q}\lbrace f_t(Q)- \langle S,Q \rangle \rbrace. 
\end{align*}
By Proposition \ref{prop:convF_grad}, it is clear that the infimum in the above equation is obtained at the values of $Q$ such that the conditions in (\ref{eq:dual_prob}b--d) are satisfied. In addition, by \eqref{subeq:Z_eqn} we further have
\begin{align*}
   \langle H, Z\rangle = \langle H,\overline{\Psi(H)} \rangle + \langle H,\hat{\Phi}_t Z\hat{\Phi}_t^{\top} \rangle,
\end{align*}
which together with \eqref{subeq:H_eqn} implies
\begin{align*}
    2f_t(Q) &= \langle H,\Psi(H) \rangle =\langle H,\overline{\Psi(H)} \rangle\\
    &=\langle H, Z\rangle-\langle H,\hat{\Phi}_t Z\hat{\Phi}_t^{\top} \rangle\\
    &=\tr (H^{\top}Z)-\tr (H^{\top}\hat{\Phi}_t Z\hat{\Phi}_t^{\top})\\
    &=\tr[(H^{\top}-\hat{\Phi}_t^{\top}H^{\top}\hat{\Phi}_t)Z]\\
    &= \tr[\diag(Q, I_{m})^{\top}Z].
\end{align*}
Since $\langle S,Q \rangle=\langle \diag(Q, 0),Z \rangle$, substituting the above expression in the definition of $\varphi_t(S)$ then completes the proof.

\subsection{Symmetric Kronecker Product}\label{appendix:svex}
    The $\svec$ operator defines an isometry between $\mathbb{S}^n$ and $\mathbb{R}^{n(n+1)/2}$.     For any $X\in\mathbb{S}^n$ with element $x_{ij}$ on its $i$th row and $j$th column ($1\leq i,j\leq n$), $\svec(X)$ is defined by 
    \begin{align*}
        [x_{11}~\sqrt{2}x_{21}~\cdots~\sqrt{2}x_{n1}~ x_{22}~\sqrt{2}x_{32}~\cdots~\sqrt{2}x_{n2}~ \cdots~x_{nn}]^{\top},
    \end{align*}
    which is associated with the vectorization operator via some matrix $U_n\in\mathbb{R}^{n(n+1)/2\times n^2}$ such that $\svec(X)=U_n\mvec(X)$ and $U_nU_n^{\top}=\mathcal{I}_n$. 
    The symmetric Kronecker product of two matrices $A,B\in\mathbb{R}^{n\times n}$ produces a square matrix of order $n(n+1)/2$, whose action on $\svec(X)$ is defined by
    $$ (A\otimes_sB)\svec(X)=\frac{1}{2}\svec(AXB^{\top}+BXA^{\top}). $$
    The following properties are used in this paper~\cite{tunccel2005strengthened}:
    \begin{enumerate}
        \item If $A$ and $B$ are positive (semi-)definite, then $A\otimes_sB$ is also positive (semi-)definite.
        \item For any $A\in \mathbb{R}^{n\times n}$, denote $\lambda_i$ as its eigenvalue and $v_i$ as the corresponding eigenvector with $i=1,2,\cdots,n$. Then $\lambda_i\lambda_j$ is an eigenvalue of $A\otimes_sA$ and $\frac{1}{2}\svec(v_iv_j^{\top}+v_jv_i^{\top})$ is the corresponding eigenvector.
    \end{enumerate}
\subsection{Proof of Proposition \ref{prop:NT_direction}}\label{append:NT_direction}
Firstly, we show that (\ref{eq:newton}a-\ref{eq:newton}c) and \eqref{eq:NT_QS} are equivalent to \eqref{eq:newton_svec}. Using the fact that $\svec(\Delta S)=U_{n}\mvec(\Delta S)$ and $\mvec(\Delta Z)=U_{n+m}^{\top}\svec(\Delta Z)$, one can directly obtain $\svec(\Delta S)=E\svec(\Delta Z)$ from \eqref{eq:newton_s} based on basic properties of the Kronecker product. Applying the symmetrized Kronecker product to both sides of \eqref{eq:newton_z} and \eqref{eq:NT_QS} yields the second and fourth rows of \eqref{eq:newton_svec} respectively. Furthermore, substituting 
\begin{equation*}
\begin{aligned}
	&\svec(\diag(\Delta Q, 0))=\svec(\begin{bmatrix}I_n & 0 \end{bmatrix}^{\top}\Delta Q \begin{bmatrix}I_n&0 \end{bmatrix})\\
=&U_{n+m}\mvec(\begin{bmatrix}I_n & 0 \end{bmatrix}^{\top}\Delta Q \begin{bmatrix}I_n&0 \end{bmatrix})=E^{\top}\svec(\Delta Q),
\end{aligned}
\end{equation*}
into \eqref{eq:newton_H} implies that $G^{\top}\svec(\Delta H)+E^{\top}\svec(\Delta Q)=0$. Hence, solving the Newton's search direction in (\ref{eq:newton_s}-\ref{eq:newton_H}) and \eqref{eq:NT_QS} is equivalent to that of \eqref{eq:newton_svec}. 

Next we show that \eqref{eq:newton_svec} admits a unique solution. It is obvious that $F\succeq 0$. Since $\hat{\Phi}_t$ is stable, all the eigenvalues of $\hat{\Phi}_t \otimes_s \hat{\Phi}_t$ are located within the unit cycle, which guarantees the invertibility of $G$. Since $W\succ 0$, we further have $W_s\succ 0$. Applying Gaussian elimination to \eqref{eq:newton_svec} yields
\begin{equation*}
    (EG^{-1}FG^{-\top}E^{\top}+W_s)\svec(\Delta Q)=\svec(\mu Q^{-1}-S),
\end{equation*}
where $\svec(\Delta Q)$ is uniquely obtained as $EG^{-1}FG^{-\top}E^{\top}+W_s\succ 0$. Then $\svec(\Delta S)$, $\svec(\Delta Z)$ and $\svec(\Delta H)$ can also be uniquely computed from the first three rows in \eqref{eq:newton_svec}.    

\subsection{Proof of Proposition \ref{prop:Phi_disturb}}\label{append:Phi_disturb}
    Recall that
\begin{equation*}
    \begin{bmatrix}
            A & B
    \end{bmatrix}=(\bar{X}_{1,t}-\bar{W}_{0,t})\Lambda_t^{-1}=(\bar{X}_{1,t+1}-\bar{W}_{0,t+1})\Lambda_{t+1}^{-1}.
\end{equation*}    
Along with \eqref{eq:WZ_bnd}, the variation of $\hat{\Phi}_{t}$ is uniformly bounded by SNR as
\begin{equation*}
\begin{aligned}
\norm{\hat{\Phi}_{t}-\Phi_{\bar{K}}}_{2} &=\norm{\begin{bmatrix}I&-\bar{K}^{\top}\end{bmatrix}^{\top}\bar{W}_{0,t}\Lambda_t^{-1}}_2\\
&\leq (\norm{\bar{K}}_2+1)\frac{2\kappa}{\gamma\zeta} =\frac{2h_{1}\kappa}{\gamma\zeta}.
\end{aligned}
\end{equation*}
By \eqref{eq:data_recur}, we further have
\begin{equation*}
    \begin{aligned}
        &\norm{\bar{X}_{1,t+1}\Lambda_{t+1}^{-1}-\bar{X}_{1,t}\Lambda_{t}^{-1}}_2\\ =& \norm{\frac{1}{t}x_{t+1}z_t^{\top}\Lambda_{t}^{-1}-\bar{X}_{1,t+1}\frac{\Lambda_{t}^{-1}z_tz_t^{\top}\Lambda_{t}^{-1}}{t+z_t^{\top}\Lambda_{t}^{-1}z_t}}_2\\
        \leq&\frac{1}{t}\norm{x_{t+1}}\norm{z_{t}}\norm{\Lambda_{t}^{-1}}_2+\norm{\bar{X}_{1,t+1}}_2\norm{\frac{\Lambda_{t}^{-1}z_tz_t^{\top}\Lambda_{t}^{-1}}{t+z_t^{\top}\Lambda_{t}^{-1}z_t}}_2.
    \end{aligned}
\end{equation*}
Since $\Lambda_t \succ0$, it holds that
\begin{align*}
    &\norm{\frac{\Lambda_{t}^{-1}z_tz_t^{\top}\Lambda_{t}^{-1}}{t+z_t^{\top}\Lambda_{t}^{-1}z_t}}_2\leq \frac{\norm{\Lambda_{t}^{-1}}_2\norm{\Lambda_{t}^{-1/2}z_tz_t^{\top}\Lambda_{t}^{-1/2}}_2}{t+z_t^{\top}\Lambda_{t}^{-1}z_t}\\
    \leq &\norm{\Lambda_{t}^{-1}}_2\frac{z_t^{\top}\Lambda_{t}^{-1}z_t}{t+z_t^{\top}\Lambda_{t}^{-1}z_t}\leq\frac{1}{t}\norm{\Lambda_{t}^{-1}}_2^2\norm{z_t}^2.
\end{align*}
Under Assumption \ref{assump:PE}, we have $\underline{\sigma}(Z_{0,t})\geq\sqrt{t}\gamma\zeta/2$ \cite{ZhaoData2025}, which implies $\overline{\sigma}(\Lambda_t^{-1})\leq{4}/{(\gamma\zeta)^2}$. Since $\lbrace x_t\rbrace_{t\geq 0}$ and $\lbrace u_t\rbrace_{t\geq 0}$ are upper bounded, it is clear that $z_t$ and $\norm{\bar{X}_{1,t+1}}_2$ are also bounded. Equipped with the above results, we can  obtain 
\begin{equation*}\label{eq:Delta_Phi_t_norm_bounds}
\begin{aligned}
\norm{\Delta_{\hat{\Phi},t}}_{2} =&\norm{\begin{bmatrix}I&-\bar{K}^{\top}\end{bmatrix}^{\top}(\bar{X}_{1,t+1}\Lambda_{t+1}^{-1}-\bar{X}_{1,t}\Lambda_t^{-1})}_2 \\
\leq& (\norm{\bar{K}}_2+1)\norm{\bar{X}_{1,t+1}\Lambda_{t+1}^{-1}-\bar{X}_{1,t}\Lambda_t^{-1}}_2 \leq \frac{h_4}{t},
\end{aligned}
\end{equation*}
for some $h_4>0$ that is determined by $\bar{K}, \gamma,\zeta$ and the magnitudes of exciting signals.

\subsection{Proof of Lemma \ref{prop:H_perturb_bnd}}\label{append:H_perturb_bnd}
Firstly, if $\frac{\gamma}{\kappa} \geq \tau_{2}$, for any $t\geq t_0$ we have $8h_3h_4(h_2+\frac{2h_1\kappa}{\gamma\zeta})\leq t$, which implies
\begin{equation}\label{eq:tau2_DeltaPhi_upperbounds}
\begin{split}
\norm{\Delta_{\hat{\Phi},t}}_{2} & \leq  \frac{h_4}{t} 
\leq \frac{1}{8h_{3}(h_{2} + h_{1}\frac{2\kappa}{\gamma \zeta})} \leq \frac{1}{4 \norm{\Sigma_{\hat{\Phi}_{t}}}_{2}(1+\norm{\hat{\Phi}_{t}}_{2})}.
\end{split} 
\end{equation}
Then by Lemma \ref{lemma:lyapunov_disturbation}, it holds that
\begin{equation*}
\begin{split}
\norm{H_{t}-\bar{H}_{t}}_{2} & \leq 4 \norm{\Sigma_{\hat{\Phi}_{t}}}_{2}^{2}(1 + \norm{\hat{\Phi}_{t}}_{2})(1 + \norm{Q_{t}}_{2}) \norm{\Delta_{\hat{\Phi},t}}_{2} \\ & \leq \frac{16}{t}h_{3}^{2}h_{4}(h_{2} + h_{1}\frac{2\kappa}{\gamma \zeta})  (1 + \norm{Q_{t}}_{2}).
\end{split}
\end{equation*}
Secondly, according to \cite{tippett1999upper}, we can obtain that
\begin{equation*}
\norm{H_{t}}_{2} \leq (\norm{Q_{t}}_{2} + 1) \norm{\Sigma_{\hat{\Phi}_{t}}}_{2} \leq 2(1 + \norm{Q_{t}}_{2}) h_{3}.
\end{equation*}

\subsection{Proof of Lemma \ref{proposition:normbound_Zt}}\label{append:normbound_Zt}
We introduce $Z_{t_{1}}$ such that
\begin{equation*}
Z_{t_{1}} =  \overline{\Psi(H_{t})}+ \hat{\Phi}_{t+1} Z_{t_{1}} \hat{\Phi}_{t+1}^{\top}.
\end{equation*}
Then $\norm{Z_{t} - \bar{Z}_{t}}_{2} \leq \norm{Z_{t} - Z_{t_{1}}}_{2} + \norm{Z_{t_{1}} - \bar{Z}_{t}}_{2}$.
Firstly, if $\frac{\gamma}{\kappa} \geq \tau_{2}$, \eqref{eq:tau2_DeltaPhi_upperbounds} is still valid. By Lemma \ref{lemma:lyapunov_disturbation}, we have
\begin{equation}\label{eq:Zt_perturbation_first}
\begin{aligned}
& \norm{Z_{t} - Z_{t_{1}}}_{2} \leq 16h_{3}^{2}(h_{2} + h_{1}\frac{2\kappa}{\gamma \zeta}) \norm{\overline{\Psi(H_{t})}}_{2} \norm{\Delta_{\hat{\Phi}}}_{2}  \\ \leq & \frac{16}{t}h_4h_{3}^{2}(h_{2} + h_{1}\frac{2\kappa}{\gamma \zeta})  h_{1}^{2} \norm{H_{t}}_{2} \\ \leq & \frac{32}{t} h_{1}^{2}h_{3}^{3}h_4(h_{2} + h_{1}\frac{2\kappa}{\gamma \zeta})  (1 + \norm{Q_{t}}_{2}). 
\end{aligned}
\end{equation}
Furthermore, the term $Z_{t_{1}} - \bar{Z}_{t}$ satisfies the equation:
\begin{equation*}
Z_{t_{1}} - \bar{Z}_{t} = \overline{\Psi(H_{t})} - \overline{\Psi(\bar{H}_{t})} + \hat{\Phi}_{t+1}(Z_{t_{1}} - \bar{Z}_{t})\hat{\Phi}_{t+1}^{\top}.
\end{equation*}
Since $\overline{\Psi(H_{t})} - \overline{\Psi(\bar{H}_{t})}$ may not be positive definite, we construct the following Lyapunov equation: 
\begin{equation*}
Z_{t_{2}} = \norm{\overline{\Psi(H_{t})} - \overline{\Psi(\bar{H}_{t})})}_2 I_{n} + \hat{\Phi}_{t+1}Z_{t_{2}}\hat{\Phi}_{t+1}^{\top}.
\end{equation*}
By monotonicity of the Lyapunov equation, it is clear that $ Z_{t_{2}} \succeq  Z_{t_{1}} - \bar{Z}_{t}$, which indicates that
\begin{equation}\label{eq:Zt_perturbation_second}
\begin{aligned}
\norm{Z_{t_{1}} - \bar{Z}_{t}}_{2} & \leq \norm{Z_{t_{2}}}_{2} \leq \norm{\overline{\Psi(H_{t})} - \overline{\Psi(\bar{H}_{t})}}_{2} \norm{\Sigma_{\hat{\Phi}_{t+1}}}_{2} \\ & \leq 2h_{3} h_{1}^{2} \norm{H_{t} - \bar{H}_{t}}_{2} \\ & \leq \frac{32}{t} h_{1}^{2}h_{3}^{3}h_4(h_{2} + h_{1}\frac{2\kappa}{\gamma \zeta})  (1 + \norm{Q_{t}}_{2}). 
\end{aligned}    
\end{equation}
Combining \eqref{eq:Zt_perturbation_first} and \eqref{eq:Zt_perturbation_second}, the conclusion follows.

\section*{References}
\vspace{-0.4cm}
\bibliographystyle{IEEEtran}
\bibliography{ref}

\vspace{-1cm}
\begin{IEEEbiography}[{\includegraphics[width=1in,height=1.25in,clip,keepaspectratio]{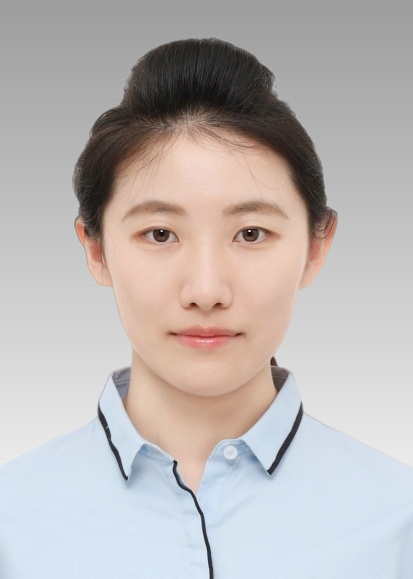}}]{Yibei Li} Yibei Li received her B.E. degree in Automation from Harbin Institute of Technology, China, in 2015. She received her Licentiate degree and Ph.D. degree in  Optimization and Systems Theory from
the Department of Mathematics, KTH Royal Institute of Technology, Sweden, in 2019 and 2022 respectively. From 2022 to 2024, she was a Wallenberg-NTU Postdoctoral Fellow at Nanyang Technological University, Singapore. She currently is an associate professor at the Academy of Mathematics and Systems Science, Chinese Academy of Sciences, China. Her research interests include inverse optimal control and estimation, nonlinear control theory, game theory and multi-agent systems. 
\end{IEEEbiography}
\vspace{-1.1cm}
\begin{IEEEbiography}[{\includegraphics[width=1in,height=1.25in,clip,keepaspectratio]{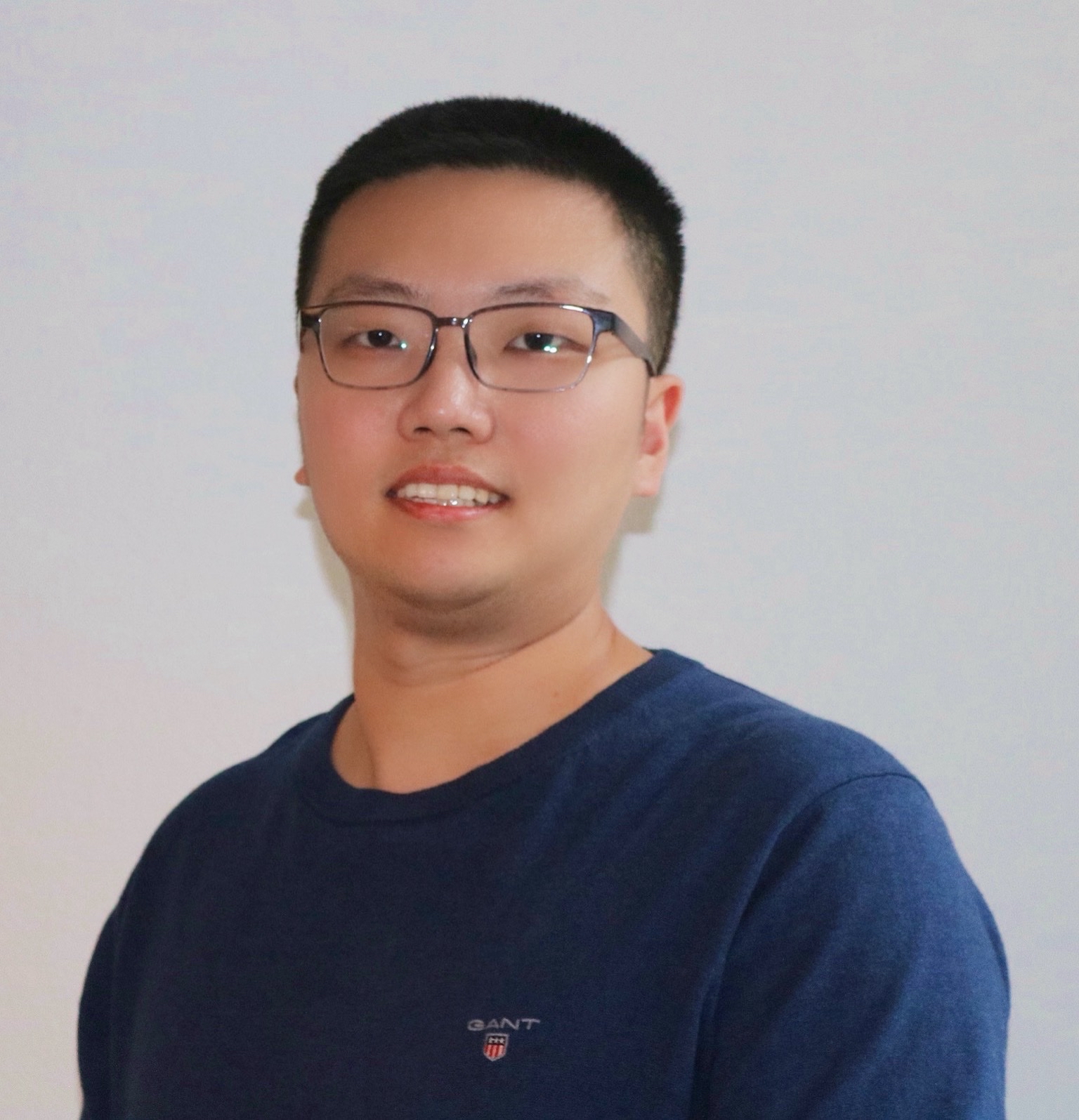}}]{Yuexin Cao} Yuexin Cao received the B.Sc. degree from Shanghai University of Finance and Economics, Shanghai, China, in 2018, and the M.Sc. degree from Fudan University, Shanghai, China, in 2020, both in Applied Statistics. He is currently pursuing the Ph.D. degree in the Department of Mathematics, KTH Royal Institute of Technology, Stockholm, Sweden. His research interests include inverse optimal control, and optimization of networked systems.
\end{IEEEbiography}
\vspace{-1.3cm}
\begin{IEEEbiography}[{\includegraphics[width=1in,height=1.25in,clip,keepaspectratio]{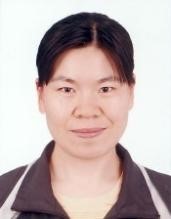}}]{Zhixin Liu}
Zhixin Liu received the B.S. degree in mathematics from Shandong University in 2002, and the Ph.D. degree in control theory from Academy of Mathematics and Systems Science (AMSS), Chinese Academy of Sciences (CAS) in 2007.  She is currently a Professor of AMSS, CAS, and the director of Key Laboratory of Systems and Control, CAS.
She had visiting positions with the KTH Royal Institute of Technology, Stockholm, Sweden, University of New South Wales, Canberra, Australia, and University of Maryland, College Park, MD, USA. Her current research interests include multi-agent systems, distributed control, and distributed estimation and filtering.
\end{IEEEbiography}
\vspace{-1.1cm}
\begin{IEEEbiography}[{\includegraphics[width=1in,height=1.25in,clip,keepaspectratio]{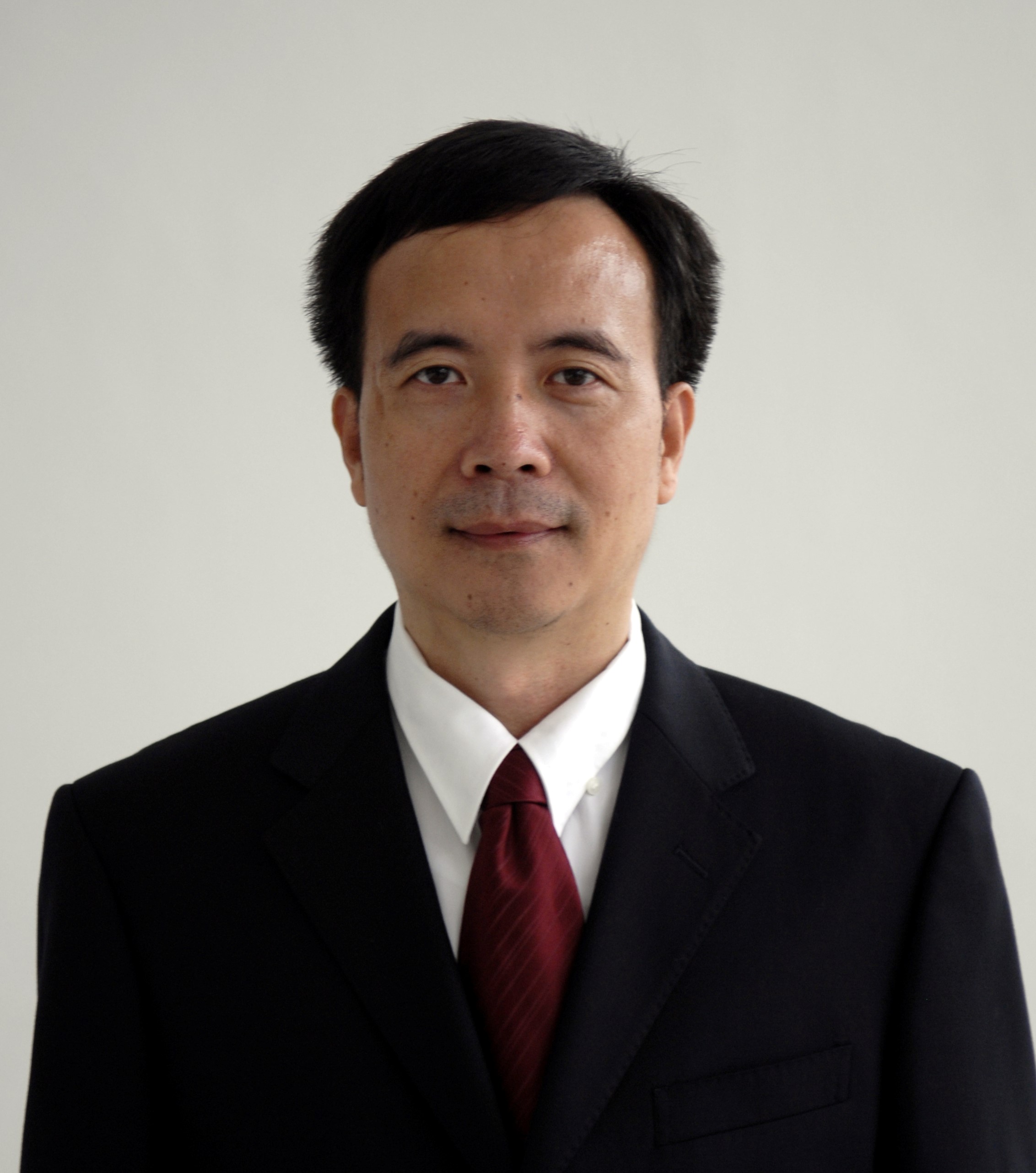}}]{Lihua Xie} (Fellow, IEEE) Lihua Xie received the Ph.D. degree in electrical engineering from the University of Newcastle, Australia, in 1992. Since 1992, he has been with the School of Electrical and Electronic Engineering, Nanyang Technological University, Singapore, where he is currently President's Chair in Control Engineering and Director, Center for Advanced Robotics Technology Innovation. He served as the Head of Division of Control and Instrumentation and Co-Director, Delta-NTU Corporate Lab for Cyber-Physical Systems. He held teaching appointments in the Department of Automatic Control, Nanjing University of Science and Technology from 1986 to 1989. 

Dr. Xie’s research interests include robust control and estimation, networked control systems, multi-agent networks, and unmanned systems. He is an Editor-in-Chief for Unmanned Systems and has served as Editor of IET Book Series in Control and Associate Editor of a number of journals including IEEE Transactions on Automatic Control, Automatica, IEEE Transactions on Control Systems Technology, IEEE Transactions on Network Control Systems, and IEEE Transactions on Circuits and Systems-II. He was an IEEE Distinguished Lecturer (Jan 2012 – Dec 2014). Dr. Xie is Fellow of Academy of Engineering Singapore, IEEE, IFAC, and CAA.
\end{IEEEbiography}

\end{document}